\documentclass[11pt,a4paper]{amsart}

\usepackage{graphicx}
\usepackage{nicefrac}
\usepackage{amsmath}
\usepackage{amsfonts}
\usepackage{amssymb}
\usepackage{amsthm}
\usepackage{pst-all}
\usepackage{fancybox}
\usepackage{url}
\usepackage{mathrsfs}
\usepackage{yhmath}
\usepackage{subfigure}
\input xy
\xyoption{all}

\newtheorem{theorem}{Theorem}
\newtheorem*{theorem*}{Theorem}

\newtheorem{definition}[theorem]{Definition}
\newtheorem{lemma}[theorem]{Lemma}
\newtheorem*{lemma*}{Lemma}

\newtheorem{proposition}[theorem]{Proposition}
\newtheorem*{proposition*}{Proposition}

\theoremstyle{definition}
\newtheorem{remark}[theorem]{Remark}
\newtheorem{example}[theorem]{Example}

\thanks{The author is funded by the Spanish Ministerio de Ciencia e Innovaci\'on through grants PGC2018-098321-B-I00 and RYC2018-025843-I.}

\begin{document}

\author{J. J. S\'anchez-Gabites}
\address{Facultad de Ciencias Matem{\'a}ticas. Universidad Complutense de Madrid. 28040 Madrid (Espa{\~{n}}a)}
\email{jajsanch@ucm.es}

%\keywords{Invariant knots and links, Low dimensional dynamics, }

%\title[]{A criterion to establish that an invariant set of a flow in $\mathbb{R}^3$ has a nontrivial homology}

\title[]{A criterion to detect a nontrivial homology of an invariant set of a flow in $\mathbb{R}^3$}

\begin{abstract} Consider a flow in $\mathbb{R}^3$ and let $K$ be the biggest invariant subset of some compact region of interest $N \subseteq \mathbb{R}^3$. The set $K$ is often not computable, but the way the flow crosses the boundary of $N$ can provide indirect information about it. For example, classical tools such as Wa\.{z}ewski's principle or the Poincaré-Hopf theorem can be used to detect whether $K$ is nonempty or contains rest points, respectively. We present a criterion that can establish whether $K$ has a nontrivial homology by looking at the subset of the boundary of $N$ along which the flow is tangent to $N$. We prove that the criterion is as sharp as possible with the information it uses as an input. We also show that it is algorithmically checkable.
\end{abstract}

\maketitle

\section{Introduction}

Let $\varphi$ be a continuous flow on $\mathbb{R}^3$. We focus on some compact region of interest $N \subseteq \mathbb{R}^3$ in phase space and want to obtain information about the largest invariant set $K$ it contains. Typically $K$ cannot be computed explicitly and indirect methods are needed. For instance, Wa\.zewski's principle (\cite{wazewski1}) can be used to analyze whether $K \neq \emptyset$, or in a smooth setting the Poincaré-Hopf theorem can detect the presence of rest points in $K$. The important feature of these tools is that they only require that one knows how the trajectories of the flow cross the boundary of $N$ and this is usually computable. For example, if the flow is generated by a vectorfield $X$ through a differential equation then $X$ alone already provides this data, without the need to integrate the equation.

The information about $K$ provided by these tools is limited. In this paper we obtain a criterion that still uses the same sort of data as an input but detects whether $K$ must have a nontrivial one-dimensional homology. Thus, for instance, it can distinguish a situation where $K$ is a rest point from another one where $K$ is a rest point together with a homoclinic orbit. Obtaining homological information about $K$ seems a natural problem and has already been studied in \cite{conleyeaston1}, \cite{easton2}, \cite{gierzkiewicz1}, \cite{pozniak1}, to cite a few. The main difficulty is, of course, that the region $N$ and its maximal invariant subset $K$ do not generally bear any direct geometric relationship whatsoever. The aforementioned papers are more ambitious than ours because they are valid in arbitrary dimensions or look for lower bounds on the rank of the homology of $K$, rather than just showing that it is nontrivial. However, for the setting considered in this paper our results are sharper, and in fact we will prove that they are the sharpest possible given the data they take as an input.
\medskip

Following the literature cited above (especially \cite{conleyeaston1}) we will focus on a certain type of regions $N$ called isolating blocks. To motivate the definition suppose for a moment that the flow is generated by a smooth vectorfield $X$ and $N$ is a smooth manifold. At every point $p \in \partial N$ the flow either enters or exits $N$ transversally or is tangent to $\partial N$. The latter happens precisely when $X(p) \cdot \nu(p) = 0$, where $\nu$ is a normal vectorfield along $\partial N$. Perhaps after perturbing $N$ (hence $\nu$) very slightly we can achieve that $X$ and $\nu$ be transverse, and then the implicit equation $X \cdot \nu = 0$ determines a finite collection of mutually disjoint simple closed curves in $\partial N$. We call these the tangency curves of $N$, or $t$--curves for short, and use the letter $\tau$ to denote them generically. The manifold $N$ is called an isolating block when all tangencies are external; i.e. whenever the flow is tangent to $\partial N$ at some $p$, a short portion of the trajectory centered at $p$ is disjoint from $N$ except at the tangency point $p$ itself. The definition in the topological case is a straightforward adaptation of this; it is given in Section \ref{sec:preliminaries} together with some remarks about how restrictive the ``no interior tangencies'' condition is.

\begin{figure}[h]
\begin{pspicture}(0,0)(11,2.3)
   % \psgrid(0,0)(11,1.8)
\rput[bl](0,0.5){\scalebox{0.6}{\includegraphics{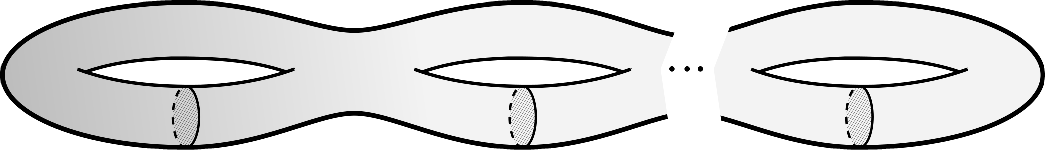}}}
\rput[tl](10.2,0.6){$H_g$}
\rput[tl](1.8,0.3){$D_1$}
\rput[tl](5.2,0.3){$D_2$}
\rput[tl](8.8,0.3){$D_g$}
\end{pspicture}
\caption{\label{fig:standard_hdbdy}}
\end{figure}

Before stating our results we recall one last definition. A handlebody $N$ is a compact $3$--manifold homeomorphic to the standard model $H_g$ of Figure \ref{fig:standard_hdbdy}. The number $g \geq 0$ is called the genus of the handlebody. A handlebody of genus $0$ is a $3$--ball; a handlebody of genus $1$ is a solid torus. Notice that $H_g$ is only an abstract model; a handlebody $N \subseteq \mathbb{R}^3$ can be knotted.

There is a more intrinsic definition of handlebodies. Let $N$ be a compact $3$--manifold with boundary. A $2$--disk $D \subseteq N$ is properly embedded if $D \cap \partial N = \partial D$. The disks $\{D_i\}$ in Figure \ref{fig:standard_hdbdy} are all properly embedded. Notice that if $H_g$ is cut along these disks the resulting manifold is a $3$--ball. This is a defining property of handlebodies: $N$ is a handlebody of genus $g$ if and only if it contains a collection of $g$ mutually disjoint, properly embedded disks $\{D_i\}$ such that when $N$ is cut along these disks one obtains a $3$--ball. The collection $\{D_i\}$ is called a (complete) cut system for $N$.
\medskip

Our first result is the following: 

\begin{theorem} \label{thm:main1} Let $N$ be an isolating block for a compact invariant set $K$. Assume that $\partial N$ is connected (i.e. that $N$ has no ``cavities''). Then, if the one-dimensional cohomology of $K$ is trivial, the following must hold:
\begin{itemize}
    \item[(i)] $N$ is a handlebody.
    \item[(ii)] There exists a complete system of cutting disks $\{D_i\}$ for $N$ such that each curve $\partial D_i$ intersects the set of tangency curves $\bigcup \tau$ transversally at two points.
\end{itemize}
\end{theorem}

Here cohomology means \v{C}ech cohomology. In this context it is usually preferred over singular homology since it is better suited to compacta having a bad local structure (as may well be the case with an invariant compact set). One should think of Theorem \ref{thm:main1} as a criterion to ensure that $K$ has a nontrivial one-dimensional cohomology: this will be the case whenever $N$ is not a handlebody, or it is a handlebody but does not have a complete cut system as described in (ii).

\begin{example} \label{ejem:basico} Suppose we observe one of the solid tori $N$ depicted in Figure \ref{fig:tori} as an isolating block. Here, and in most subsequent figures, we follow the convention of painting transverse entry points white (or light gray) and transverse exit points dark gray.

For a solid torus a complete cut system consists of a single meridional disk which is unique up to isotopy, and this makes it easy to check condition (ii) in Theorem \ref{thm:main1}. In panel (a) any meridian of $\partial N$, however contorted, must intersect the four parallel tangency curves. Thus by Theorem \ref{thm:main1} the maximal invariant subset of $N$ must have a nontrivial one-dimensional cohomology. The solid torus $N$ in panel (b) admits a slightly bent meridional disk $D$ whose boundary (shown in a dotted outline) does satisfy condition (ii) of Theorem \ref{thm:main1}, so we cannot draw any conclusion about the cohomology of its maximal invariant subset.

\begin{figure}[h]
\null\hfill
\subfigure[]{
\scalebox{0.9}{\begin{pspicture}(0,0)(5.8,2.2)
	%\psgrid(0,0)(5.8,2.2)
	\rput[bl](0,0){\scalebox{0.5}{\includegraphics{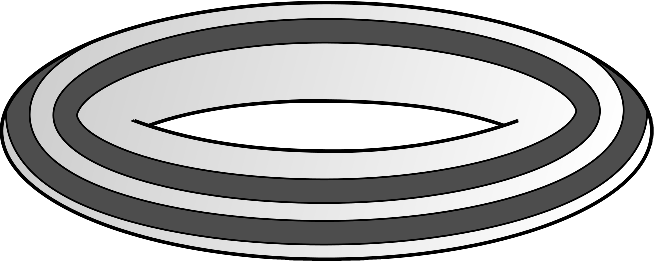}}}
\end{pspicture}}}
\hfill
\subfigure[]{
\scalebox{0.9}{\begin{pspicture}(0,0)(5.8,2.2)
	%\psgrid(0,0)(5.8,2.2)
	\rput[bl](0,0){\scalebox{0.5}{\includegraphics{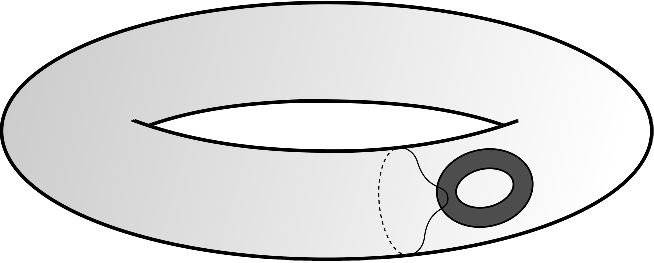}}}
	\rput[r](3.1,0.5){$\partial D$}
\end{pspicture}}}
\hfill\null
\caption{\label{fig:tori}}
\end{figure}

 The example in panel (a) is the simplest one where the condition of Theorem \ref{thm:main1}.(ii) is not satisfied. It was already considered by Conley and Easton (\cite[Section 3, pp. 47ff.]{conleyeaston1}). We have constructed the example in panel (b) so that both have the same Conley index (namely $\mathbb{S}^1 \vee \mathbb{S}^2$); this implies that other classical tools such as the Poincaré-Hopf theorem cannot tell apart the two situations (\cite{mccord3}).
\end{example}

Recall that a closed curve in $N$ is called essential if it is not contractible. One has the following easy generalization of Example \ref{ejem:basico}.(a): 

\begin{example} \label{ex:number} Let an isolating block $N$ be a handlebody of genus $g$. If the number of essential tangency curves in $N$ is $\geq 2g$, then $K$ has a nontrivial cohomology.

To prove this, let $\{D_i\}$ be a complete cut system for $N$ and notice that if a tangency curve $\tau$ does not intersect any of the $\partial D_i$ then after cutting $N$ open along the $D_i$ the curve $\tau$ would remain, untouched, in the boundary of resulting $3$--ball. Since it is evidently contractible in that ball, it is also contractible in $N$. It follows that every essential tangency curve in $N$ must intersect at least one of the $\partial D_i$. Thus if there are more than $2g$ of those curves in number, one of the $\partial D_i$ must intersect the system of tangency curves more than twice and the result follows from Theorem \ref{thm:main1}.(ii). 
\end{example}

 It is interesting to observe that, at least in the smooth setting, the nontrivial cohomology detected by Theorem \ref{thm:main1} is stable under small perturbations of the flow:

\begin{remark} Let $N$ be an isolating block for a smooth flow $\varphi$. Suppose that by using the criterion provided by (the counterpositive of) Theorem \ref{thm:main1} we know that the maximal invariant subset $K$ of $N$ has a nontrivial one-dimensional cohomology. Then there exists $\delta > 0$ such that for any smooth flow $\varphi'$ which is $\delta$-close to $\varphi$ over $N$, the maximal invariant subset of $N$ (with respect to $\varphi'$) also has a nontrivial one-dimensional cohomology.    
\end{remark}

The condition that $\varphi$ and $\varphi'$ be $\delta$-close over $N$ means that $\|X(p) - X'(p)\| < \delta$ for every $p \in N$, where $X$ and $X'$ are the velocity vectorfields of the flows. With $N$ and $\varphi$ as above, a stability theorem of Conley and Easton (\cite[Theorem 1.6, p. 39]{conleyeaston1}) states the following: for every $\epsilon > 0$ there exists $\delta > 0$ with the property that for any flow $\varphi'$ which is $\delta$-close to $\varphi$ over $N$ there exists a homeomorphism $h$ of $\mathbb{R}^3$ that moves points less than $\epsilon$ and carries $N$ onto an isolating block $N' := h(N)$ for $\varphi'$; moreover, $h$ carries the tangency curves of $N$ onto the tangency curves of $N'$. Since both conditions of Theorem \ref{thm:main1} are invariant under such a homeomorphism and $N$ failed to satisfy at least one of them by assumption (for $\varphi$), the same is true of $N'$ (for $\varphi'$). Thus the maximal invariant subset of $N'$ for $\varphi'$ must have a nontrivial one-dimensional cohomology. By choosing $\epsilon$ sufficiently small (so that $N$ and $N'$ are very close to each other) and perhaps reducing $\delta$ we can ensure that the maximal invariant subset of $N$ and of $N'$ for $\varphi'$ is the same, and the remark follows.
\medskip

Next we consider how sharp the criterion provided by Theorem \ref{thm:main1} is. To discuss this it is convenient to introduce a definition. Let $N$ be a compact $3$--manifold. A colouring of $N$ means a decomposition $\partial N = P \cup Q$ where $P,Q$ are compact $2$--manifolds (possibly empty) with a common boundary $\partial P = \partial Q = P \cap Q$. We think of the interiors of $P$ and $Q$ as being painted white and gray, say. Clearly $P \cap Q$ is a finite collection of disjoint, simple closed curves $\tau$ which we call the $t$--curves of $N$. Of course, this is just a topological abstraction of the ``entry, exit, tangency'' information carried by the boundary of an isolating block. Given a coloured manifold $(N,P,Q)$ in $\mathbb{R}^3$, we say that a flow in $\mathbb{R}^3$ realizes $N$ as an isolating block if every $p$ lying on a $t$--curve is an exterior tangency and every point in the interior of $P$ (resp. of $Q$) is a transverse entry (resp. exit) point.

\begin{theorem}\label{thm:main2} Let $(N,P,Q) \subseteq \mathbb{R}^3$ be a (tame) coloured manifold that satisfies conditions (i) and (ii) in Theorem \ref{thm:main1}. Then there exists a flow on $\mathbb{R}^3$ which realizes $N$ as an isolating block whose maximal invariant subset $K$ is a single rest point.
\end{theorem}

This means that the criterion provided by Theorem \ref{thm:main1} is as sharp as the information about the colouring of an isolating block $N$ allows. For instance, in Example \ref{ejem:basico} we did not reach any conclusion about the isolating block $N$ shown in panel (b); now we know that there actually exists a flow in $\mathbb{R}^3$ which realizes $N$ as an isolating block for a single rest point. Notice also the following. If we start with an isolating block $N$ (with a connected boundary) whose maximal invariant subset $K$ has a trivial one-dimensional cohomology, by applying Theorems \ref{thm:main1} and \ref{thm:main2} succesively we conclude that there exsits a flow in $\mathbb{R}^3$ which realizes $N$ as an isolating block with the same colouring (i.e. with the same set of tangency curves and transverse entry and exit sets) but whose maximal invariant subset is a single rest point. In other words, there is no way to distinguish a maximal invariant subset $K$ with a trivial one-dimensional cohomology from a single rest point only using information from the boundary of an isolating block.
\medskip

Of the two conditions in Theorem \ref{thm:main1} it is (ii) that is difficult to check in practice because cut systems of a handlebody $N$ are highly non-unique (for genus $g \geq 2$). We shall call this condition the ``geometric criterion'' for brevity. When $N$ is a $3$--ball it evidently satisfies the condition regardless of its colouring, since complete cut systems are empty. When $N$ is a solid torus one can easily show that it satisfies the geometric criterion if and only if it contains at most two essential $t$--curves. For handlebodies of genus $g \geq 2$ a hands-on analysis becomes almost impossible. Our last theorem provides an algorithm that checks whether a coloured handlebody satisfies the geometric criterion. We need a preliminary discussion.

Let $N$ be a handlebody with an entirely arbitrary complete cut system $\{D_1,\ldots,D_g\}$. Assign to each disk $D_i$ a little arrow transverse to it which will indicate a positive crossing direction. Now, given any oriented simple closed curve $s$ in $N$, we can manufacture a word $W$ in the letters $x_1^{\pm 1},\ldots,x_g^{\pm 1}$ by travelling once along the curve and recording all encounters with the cutting disks by writing $x_i^{\pm 1}$ whenever $s$ crosses $D_i$, with the exponent $\pm 1$ indicating whether the crossing takes place in the positive direction. If the curve $s$ is disjoint from all the disks in the cut system the word $W$ is empty and we write $W = 1$. (The reader might recognize this as a procedure to express the free homotopy class of $s$ in the fundamental group of $N$, which is a free group on the generators $\{x_i\}$).

Given any (finite) set $S$ of words in the letters $x_i^{\pm 1}$ there is a certain purely combinatorial algorithm, called Whitehead reduction, that returns another set of words $S_{\rm min}$ in the same letters. Very roughly speaking, this process attempts to reduce the length of the words in $S$ by performing certain substitutions and cancellations; the set $S_{\rm min}$ it returns has the property that the total length of its words cannot be reduced further. We will describe this in Section \ref{sec:proof3}. For the moment this rough idea is enough to state our last theorem:

\begin{theorem} \label{thm:main3} Let $N$ be a coloured handlebody having at least one $t$--curve. Orient its $t$--curves as the boundary of the dark gray region (say). Choose any complete cut system $\{D_i\}$ and read the $t$--curves as explained above to produce a set of words $S$ in the letters $x_i^{\pm 1}$. Perform Whitehead reduction to obtain a new set $S_{\rm min}$. Then $N$ satisfies the geometric criterion (i.e. condition (ii) of Theorem \ref{thm:main1}) if and only if $S_{\rm min}$ satisfies the following:

$(A)$ For every $i$, the letters $x_i^{\pm 1}$ either do not appear at all among the words in $S_{\rm min}$ or both appear, exactly once each.
\end{theorem}

The advantage of the algebraic criterion $(A)$ over the geometric criterion of Theorem \ref{thm:main1} is that it is no longer formulated in existential terms and is algorithmically checkable. Notice also that the assumption that $N$ has at least one $t$--curve is inconsequential since otherwise $N$ cannot satisfy the geometric criterion anyway (unless it is a $3$-ball). In fact there is a good dynamical explanation for this. If $N$ has no $t$--curves it must be painted in just one colour; say white. This means that it will be positively invariant for any flow $\varphi$ that realizes it, and then its maximal invariant subset $K$ is an attractor with $H^*(K) = H^*(N) \neq 0$ (a result of Hastings \cite{hastings1}) regardless of the details of the flow.

We illustrate Theorem \ref{thm:main3} with an example:

\begin{example} Suppose we observe an isolating block $N$ which is a handlebody of genus $2$ with the system of tangency curves shown (in a top view of the handlebody) in Figure \ref{fig:ejemplo_complicado}.

\begin{figure}[h]
\begin{pspicture}(0,0)(7.5,5)
	%\psgrid(0,0)(7.5,5)
	\rput[bl](0.5,0.5){\scalebox{0.4}{\includegraphics{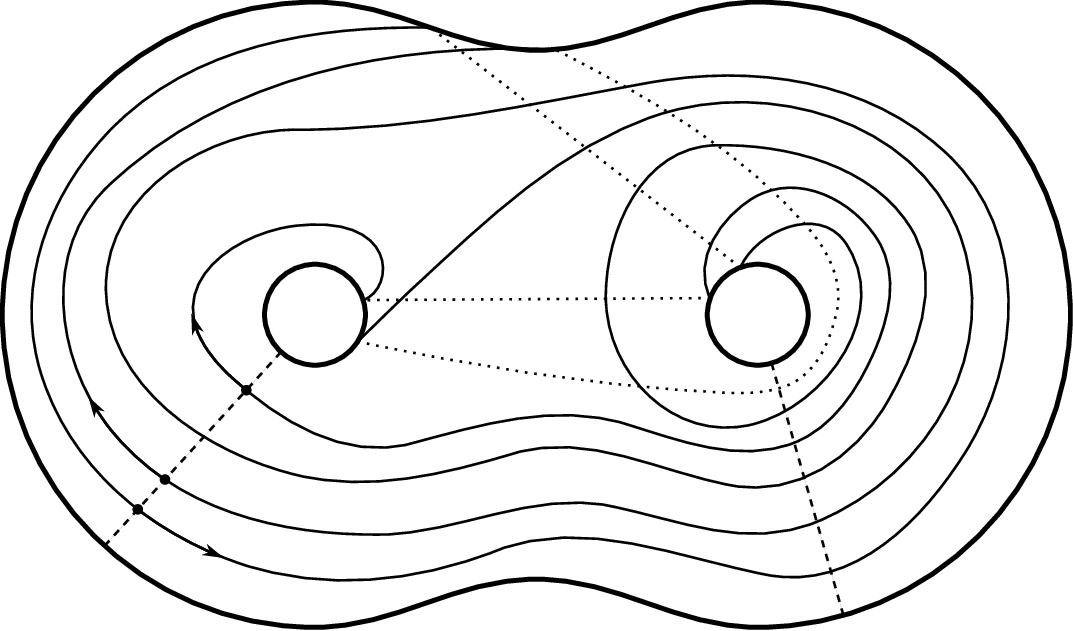}}}
         \rput[bl](0.8,0.6){$D_1$}
        \rput[bl](6.2,0.2){$D_2$}
\end{pspicture}
\caption{\label{fig:ejemplo_complicado}}
\end{figure}

We choose the standard cut system $\{D_1,D_2\}$ represented by the dashed radial lines. It is clear that the example does not satisfy condition (ii) of Theorem \ref{thm:main1} for that particular cut system; $D_1$ and $D_2$ intersect $\bigcup \tau$ a total of $4$ and $6$ times, respectively. To apply Theorem \ref{thm:main3} we orient all three $t$--curves as indicated by the arrows; we also let counterclockwise crossings of the $D_i$ be recorded with a $+1$ exponent. We start reading each curve from the thick dot lying on $D_1$. The words spelled by the $t$--curves are then $S = \{x_1 x_2 x_1 x_2 x_2, x_1^{-1}x_2^{-1}x_2^{-1},x_1^{-1}x_2^{-1}\}$. After performing Whitehead reduction we get $S_{\rm min} = \{x_1x_2,x_1^{-1},x_2^{-1}\}$, which satisfies $(A)$. Hence condition (ii) of Theorem \ref{thm:main1} is verified for some cut system $\{D'_i\}$, although it is not at all obvious in the figure what that system might look like. Dynamically, we cannot conclude anything about the cohomology of the maximal invariant subset of $N$, since by Theorem \ref{thm:main2} there even exists a flow in $\mathbb{R}^3$ which realizes this pattern of tangency curves on $N$ and has a single rest point inside $N$ as its maximal invariant subset. (Visualizing such a flow seems difficult).
\end{example}

The paper is organized as follows. Section \ref{sec:preliminaries} contains some preliminary definitions; not many are needed. Sections \ref{sec:proof1} and \ref{sec:proof2} are devoted to proving Theorems \ref{thm:main1} and \ref{thm:main2} above, respectively. Theorem \ref{thm:main3} is proved in Section \ref{sec:proof3} after recalling some algebraic preliminaries. Finally, a brief Section \ref{sec:models} discusses ``how likely'' it is for a coloured handlebody to satisfy the geometric criterion in Theorem \ref{thm:main1} or, dynamically, for an isolating handlebody to actually isolate an invariant set with a trivial cohomology.

\section{Preliminary definitions} \label{sec:preliminaries}

\subsection{} We will denote the boundary of a manifold $N$ by $\partial N$ and its interior by ${\rm Int}\ N$. The word ``interior'' will almost always be used with this meaning and not as the topological interior of a subset of some other set. 
A manifold $N \subseteq \mathbb{R}^3$ is tame if there exists a homeomorphism of $\mathbb{R}^3$ which sends $N$ onto a polyhedral set (equivalently, onto a smooth set). Tameness is a technical condition which is, in practice, always warranted. Saying that $N$ is tame is equivalent to requiring that it can be collared in $\mathbb{R}^3$.

Since an invariant compactum $K$ may have bad local topological properties it is convenient to focus on its \v{C}ech cohomology rather than its singular homology, because the former generally captures more information. It will suffice, and actually be convenient for some proofs, to take coefficients in $\mathbb{Z}_2$. Thus from now on $H^*$ denotes \v{C}ech cohomology with $\mathbb{Z}_2$ coefficients. Notice that for a locally contractible space (such as a manifold, or a polyhedron) this coincides with singular cohomology. Details about \v{C}ech cohomology can be found in \cite{eilenbergsteenrod} or \cite{spanier1}, but the reader unfamiliar with it may just prefer to interpret $H^*$ as singular cohomology; the proofs should still make sense (mostly).

\subsection{} Let $\varphi : \mathbb{R}^3 \times \mathbb{R} \longrightarrow \mathbb{R}^3$ be a flow in $\mathbb{R}^3$. We abbreviate $\varphi(p,t)$ by $p \cdot t$. Let $N \subseteq \mathbb{R}^3$ be a compact $3$--manifold.

We shall say that a point $p \in \partial N$ is a transverse (i) entry or (ii) exit point if there exists an $\epsilon > 0$ such that either (i) $p \cdot (-\epsilon,0) \cap N = \emptyset$ and $p \cdot (0,\epsilon) \subseteq {\rm Int}\ N$, or (ii) $p \cdot (-\epsilon,0) \subseteq {\rm Int}\ N$ and $p \cdot (0,\epsilon) \cap N = \emptyset$. Similarly, $p$ is an exterior tangency if $p \cdot (-\epsilon,0) \cap N = \emptyset = p \cdot (0,\epsilon) \cap N$.

\begin{definition} \label{defn:iblocks} An isolating block $N$ is a tame compact $3$--manifold $N \subseteq \mathbb{R}^3$ whose boundary $\partial N$ is the union of two compact $2$--manifolds $N^i$ and $N^o$ (one may be possibly empty) with a common boundary $\partial N^i = \partial N^o = N^i \cap N^o$ and such that:
\begin{itemize}
	\item[(i)] every $p \in {\rm Int}\ N^i$ is a transverse entry point,
	\item[(ii)] every $p \in {\rm Int}\ N^o$ is a transverse exit point,
	\item[(iii)] every $p \in N^i \cap N^o$ is an exterior tangency.
\end{itemize}
\end{definition}

The superscripts in $N^i$ and $N^o$ stand for ``in'' and ``out'' and suggest where the flow is heading. As mentioned before, each component of $N^i \cap N^o$ is called a tangency curve, or a $t$--curve for short. Notice that the maximal invariant subset $K$ of $N$ is contained in its interior because there are no interior tangencies to $\partial N$.

\begin{remark} We make a brief comment about how restrictive the definition of an isolating block is. Suppose $K$ is a compact invariant set for a continuous flow in a $3$--manifold without boundary. Assume that $K$ is isolated in the sense of Conley; i.e. it has a compact neighbourhood $N_1$ such that $K$ is the largest invariant subset of $N_1$. Then $N_1$ contains an isolating block $N$ for $K$ (proved in \cite{mio3} by blending ideas from \cite{conleyeaston1} and \cite{gierzkiewicz1}), and so from a theoretical point of view every isolated invariant compact set in $\mathbb{R}^3$ can be analyzed with the results in this paper. Under favourable circumstances it is even possible to compute an isolating block $N$ explicitly (see \cite{mischaikow1}).
\end{remark}

For any $p \in N$ we define its exit time as \[t^o(p) := \sup\ \{t \in [0,+\infty) : p \cdot [0,t] \subseteq N\}.\] Notice that $t^o(p) = +\infty$ precisely when the forward orbit of $p$ is entirely contained in $N$; otherwise $t^o(p) < + \infty$ is the time it takes the forward orbit of $p$ to leave $N$. Because the trajectories cannot slide along $\partial N$ in an isolating block, one can easily show that $p \cdot (0,t^o(p)) \subseteq {\rm Int}\ N$ and, when $0 < t^o(p) < +\infty$, the point $p \cdot t^o(p)$ is a transverse exit point. In particular points close to $p$ will have a similar exit time; it then follows that $t^o : N \longrightarrow [0,+\infty]$ is a continuous map.

\subsection{} \label{subsec:cutting} We adapt the usual notion of ``cutting a manifold along a disk'' to coloured manifolds as follows. Let $(N,P,Q)$ be a coloured $3$-manifold with $t$--curves $\{\tau\}$. A cutting disk for $N$ means a disk $D \subseteq N$ such that:
\begin{itemize}
    \item [(i)] $D$ has a (topological) regular neighbourhood in $N$; i.e. a neighbourhood $U$ such that $(U,U \cap \partial N) \cong (D,\partial D) \times [-2,2]$ via a homeomorphism that carries $x \in D \subseteq U$ to $(x,0) \in D \times \{0\}$.
    \item [(ii)] Condition (i) implies that $D \cap \partial N = \partial D$ (i.e. $D$ is properly embedded in $N$); we further require that $\partial D$ intersects $\bigcup \tau$ exactly twice, transversally.
\end{itemize}

\begin{figure}[h!]
\null\hfill
\subfigure[]{
\begin{pspicture}(0,0)(4.8,4.8)
	%\psgrid(0,0)(4.8,4.8)
	\rput[bl](0,0){\scalebox{0.6}{\includegraphics{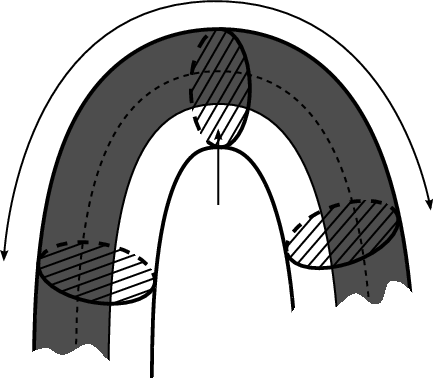}}}
        \rput[t](3.8,0.5){$N$}
	\rput[t](2.2,1.6){$D$} \rput[t](2.2,4.4){$U \cong D \times [-2,2]$}
\end{pspicture}}
\hfill
\subfigure[]{
\begin{pspicture}(0,0)(4,4)
	%\psgrid(0,0)(4,4)
	\rput[bl](0,0){\scalebox{0.6}{\includegraphics{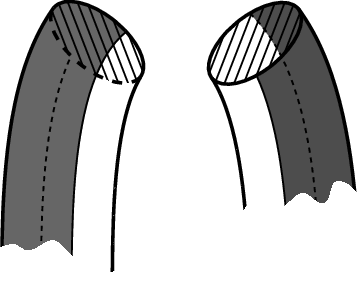}}}
        \rput[t](3.6,0.7){$N'$}
        \rput[b](0.8,3){$D_-$}
        \rput[b](3,3){$D_+$}

\end{pspicture}}
\hfill\null
\caption{\label{fig:cap}}
\end{figure}

We will often abuse notation and write $D \times [-2,2]$ instead of $U$. Condition (i) is needed to avoid wild disks, which are known to exist in $3$ dimensions. It is automatically satisfied when working in the piecewise linear or smooth category letting $U$ be a regular or tubular neighbourhood of $D$, respectively. In (ii) ``transversally'' should be understood in the obvious topological sense: if $p$ is a point of intersection of $\partial D$ and some $t$--curve $\tau$, a short arc of $\partial D$ centered at $p$ should have one of its halves on the gray side of $\tau$ and the other on its white side. By reparameterizing $U$ We can (and will) assume that for each $t \in [-2,2]$ the boundary of the slice $D \times \{t\}$ of $U$ intersects $\bigcup \tau$ in the same way as the boundary of $D$ (i.e. it intersects the same $t$--curves, and transversally). The picture conveyed by this definition is that of Figure \ref{fig:cap}.(a).

Cutting $N$ along $D$ means removing from $N$ the topological interior of $D \times (-1,1)$ to obtain a new $3$--manifold $N'$. One has \[\partial N' = \partial N \setminus (\partial D) \times [-1,1] \cup D \times \{\pm 1\}\] so most of the points of the boundary of $N'$ inherit a colouring from $N$; only points in the interior of the two ``cutting disks'' $D \times \{\pm 1\}$ remain uncoloured. We shall denote these disks by $D_{\pm}$ for brevity and often represent them in drawings filled with diagonal lines. See Figure \ref{fig:cap}.(b).

We will usually start with $N$ being a tame subset of $\mathbb{R}^3$ (or some other ambient $3$--manifold). Removing only $D \times (-1,1)$ instead of $D \times (-2,2)$ from $N$ ensures that $N'$ is still tame in $\mathbb{R}^3$.

\section{Proof of Theorem \ref{thm:main1} \label{sec:proof1}}

The basic structure of the proof goes as follows. Let $N \subseteq \mathbb{R}^3$ be an isolating block with a connected boundary. Consider the inclusion induced homomorphism $H^1(N) \longrightarrow H^1(K)$. We will show how to use (certain) elements in the kernel of this map to construct cutting disks $D$ for $N$ that are transverse to $\bigcup \tau$. Cutting $N$ along these disks produces an $N'$ such that $H^1(N') \longrightarrow H^1(K)$ is injective (i.e. all elements in its kernel have been removed). Thus if $K$ has a trivial one-dimensional cohomology so does this $N'$. Theorem \ref{thm:main1} then follows easily.

This technique was already used in \cite{mio3} to establish a certain ``handle decomposition theorem'' of isolating blocks. For the purposes of that paper we could assume that $K$ had finitely many connected components. Since here the basic tenet is that $K$ is completely unknown, this assumption has to be removed and makes the analysis slightly more delicate.

\subsection{} \label{subsec:prep} First we explain how the cutting disks will be constructed. Denote by $n^+$ the subset of points in $N^i$ whose forward trajectory never exits $N$, or equivalently $t^o(p) = +\infty$. Let $\gamma$ be an arc properly embedded in $N^i$; this means that $\gamma \cap \partial N^i = \partial \gamma$. Assume that $\gamma$ does not intersect $n^+$, so that each $p \in \gamma$ exits $N$ in a finite time $t^o(p)$. Then the set \[D := \bigcup_{p \in \gamma} p \cdot [0,t^o(p)]\] obtained by following each point in $\gamma$ until it first exits $N$ is a $2$--disk properly embedded in $N$. We say that $D$ is generated from $\gamma$ by the flow. Notice that $\partial D$ consists of the two arcs $\gamma$ and $\{p \cdot t^o(p) : p \in \gamma\}$, which intersect exactly at their endpoints. It is therefore a simple closed curve that intersects $\bigcup \tau$ exactly twice (at the endpoints of $\gamma$) and transversally.

To cut $N$ along $D$ we proceed as follows. Pick a thin strip $E \subseteq N^i$ along $\gamma$, still disjoint from $n^+$. This $E$ is a rectangle which intersects $\partial N^i$ exactly along its short sides and whose long sides run parallel to $\gamma$. Formally, $(E,E \cap \partial N^i) \cong (\gamma,\partial \gamma) \times [-1,1]$ where $\gamma$ corresponds to $\gamma \times \{0\}$. We then use the flow again to generate \[U := \bigcup_{p \in E} p \cdot [0,t^o(p)],\] which is a thickening of $D$ inside $N$ as described in Subsection \ref{subsec:cutting}, and cut $N$ along $D$ to obtain the manifold $N' := N \setminus {\rm int}_N U$. The cutting disks $D_{\pm} \subseteq \partial N'$ correspond to the disks generated by the long sides of the rectangle $E$. Notice that $N'$ contains $K$ in its interior and its entry and exit time maps are continuous, since they are restrictions of those of $N$. However $N'$ is not an isolating block because the flow slides along the cutting disks $D_{\pm}$ contained in its boundary. Every other point in $\partial N'$ retains the nature (i.e. exterior tangency, transverse entry or exit point) it had prior to the cutting.

\subsection{} Let $D \subseteq N$ be a disk generated from an arc as described above. This can be used to define a cohomology class $w_D \in H^1(N)$ which acts on a homology class $[c] \in H_1(N)$ by taking a representative $1$--chain $c$ which is transverse to $D$ and counting (modulo $2$) the number of intersections of the $1$--cycle $c$ and $D$. Formally, the $2$--cycle $D$ represents a homology class $[D] \in H_2(N,\partial N)$ and $w_D$ is its image under the Lefschetz duality isomorphism $H_2(N,\partial N) \cong H^1(N)$. It is clear that if $O$ is a neighbourhood of $K$ small enough that it is disjoint from $D$ then $w_D$ restricts to zero in $H^1(O)$; thus each $w_D$ belongs to the kernel of the inclusion induced map $H^1(N) \longrightarrow H^1(K)$. The following lemma provides a converse:

\begin{lemma} \label{lem:kernel} Assume $H^2(N,K) = 0$. Let $\gamma$ run over all arcs properly embedded in $N^i - n^+$ and let $D$ run over all the disks generated from these arcs by the flow. Then $\{w_D\}$ generates the kernel of $H^1(N) \longrightarrow H^1(K)$.
\end{lemma}

We will prove the lemma in Subsection \ref{subsec:coh_generate}. Let us admit it and continue our argument.

Focus on any one $w_D$ that is nonzero. Perform the cut-along $D$ procedure described above to obtain a new $N'$. Observe that the cohomology of $(N,N')$ is (by excision) that of $U$ relative to the cutting two disks $D_{\pm}$; i.e. that of $(D \times [-1,1], D \times \{\pm 1\})$. This cohomology is zero in degree $\neq 1$, and in degree $1$ it is $\mathbb{Z}_2$ generated by the class $\hat{w}_D \in H^1(N,N')$ which counts intersections with the disk $D$. We then have a commutative diagram  \[\xymatrix{0 & H^1(N') \ar[l] \ar[dr] & H^1(N) \ar[d] \ar[l]_-{\alpha} & H^1(N,N') = \mathbb{Z}_2(\hat{w}_D) \ar[l]_-{\beta} \\ & & H^1(K) &}\] where the row comes from the long exact sequence of the pair $(N,N')$ and the unlabeled arrows are induced by the inclusions. Since $\alpha$ is surjective, the commutativity of the triangle involving $H^1(K)$ shows that the kernel of $H^1(N') \longrightarrow H^1(K)$ is the image under $\alpha$ of the kernel of $H^1(N) \longrightarrow H^1(K)$. Notice that $\beta(\hat{w}_{D}) = w_D$ because both cohomology classes count intersections with $D$ and, since $w_D$ is nonzero by assumption, it follows that \[\dim \ker (H^1(N') \longrightarrow H^1(K)) = \dim \ker (H^1(N) \longrightarrow H^1(K)) -1.\] For later reference observe also that looking at the rest of the exact sequence for the pair $(N,N')$ and using that $\beta$ is injective we have that the inclusion $N' \subseteq N$ induces isomorphisms in $H^2$ and $H^0$.

A slight elaboration on this procedure yields the following:

\begin{proposition} \label{prop:cutsystem} Assume $H^2(N,K) = 0$. Then there exists a finite family of disjoint, properly embedded arcs $\{\gamma_i\}$ in $N^i$ such that $N$ cut along the disks $\{D_i\}$ generated by these arcs produces an $N'$ with $H^1(N') \longrightarrow H^1(K)$ injective.
\end{proposition}
\begin{proof} Assume for the sake of argument that the kernel of $H^1(N) \longrightarrow H^1(K)$ has dimension $2$. Then there exist two arcs $\gamma_1,\gamma_2$ properly embedded in $N^i$ which generate disks $D_1, D_2$ such that $\{w_{D_1},w_{D_2}\}$ generate that kernel.

{\it Case 1.} If the $\gamma_i$ are disjoint so are the disks $D_i$, and then a straightforward variation on the procedure described above shows that cutting $N$ along the $\{D_i\}$ simultaneously produces an $N'$ with the required property.

{\it Case 2.} Suppose, then, that the arcs $\gamma_i$ are not disjoint. By perturbing them slightly (which does not change the cohomology classes $w_{D_i}$) we may assume that the arcs intersect transversally at some points in the interior of $N^i$. For instance, imagine that they intersect at a single point as in Figure \ref{fig:intersect}.(a). We modify $\gamma_2$ by breaking it at the intersection point and adding arcs parallel to $\gamma_1$ to connect the two portions of $\gamma_2$ to the boundary of $N^i$ obtaining two properly embedded arcs $\gamma_2'$ and $\gamma_2''$. See Figure \ref{fig:intersect}.(b). Clearly $\gamma_2$ is homologous to $\gamma_2' + \gamma_2''$. (If there are more intersection points we do this for each of them and obtain more summands, but the argument is the same). The three arcs $\{\gamma_1,\gamma_2',\gamma_2''\}$ are now disjoint and generate classes $w_{D_1}, w_{D_2'}, w_{D_2''}$ which satisfy $w_{D_2} = w_{D_2'} + w_{D_2''}$. It follows that either $\{w_{D_1},w_{D_2'}\}$ or $\{w_{D_1},w_{D_2'}\}$ (or perhaps both) generate the kernel of $H^1(N) \longrightarrow H^1(K)$; since they are generated by disjoint arcs, we fall back onto Case 1.

\begin{figure}[h!]
\null\hfill
\subfigure[]{
\begin{pspicture}(0,0)(4,3)
	%\psgrid(0,0)(4,3)
	\rput[bl](0,0){\scalebox{0.7}{\includegraphics{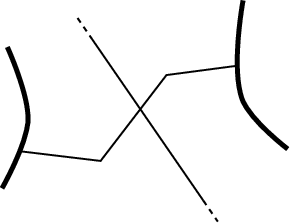}}}
	\rput[bl](3.3,1.2){$\partial N^i$}
        \rput[t](1,0.6){$\gamma_1$}
        \rput[bl](2.3,0.5){$\gamma_2$}
\end{pspicture}}
\hfill
\subfigure[]{
\begin{pspicture}(0,0)(4,3)
	%\psgrid(0,0)(4,3)
	\rput[bl](0,0){\scalebox{0.7}{\includegraphics{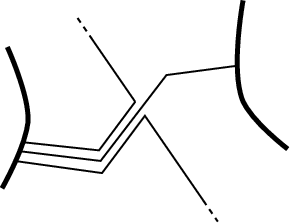}}}
 	\rput[bl](3.3,1.2){$\partial N^i$}
        \rput[bl](2.3,0.5){$\gamma_2'$}
        \rput[bl](1.3,2){$\gamma_2''$}
\end{pspicture}}
\hfill\null
\caption{\label{fig:intersect}}
\end{figure}

\end{proof}

We can now prove Theorem \ref{thm:main1} in the Introduction, which we restate here:

\begin{theorem*} Let $N \subseteq \mathbb{R}^3$ be an isolating block with a connected boundary. Assume its maximal invariant subset $K$ has $H^1(K) = 0$. Then $N$ is a handlebody and it has a complete system of cutting disks $\{D_i\}$ such that each $\partial D_i$ intersects the system of tangency curves exactly twice, and transversally.
\end{theorem*}
\begin{proof} Since $N \subseteq \mathbb{R}^3$, saying that it has a connected boundary is equivalent to saying that $H_2(N) = 0$ (write $\mathbb{S}^3 = N \cup \overline{\mathbb{S}^3 \setminus N}$ and look at the Mayer-Vietoris sequence of this decomposition) and in turn $H^2(N) = 0$ by the universal coefficient theorem. The long exact sequence in \v{C}ech cohomology for the pair $(N,K)$ and the assumption $H^1(K) = 0$ then imply $H^2(N,K) = 0$. Proposition \ref{prop:cutsystem} above yields a system of cutting disks $\{D_i\}$ such that $N$ cut along these $\{D_i\}$ produces an $N'$ with $H^1(N') \longrightarrow H^1(K)$ injective. Since $H^1(K) = 0$ by assumption, this forces $H^1(N') = 0$. We observed just before Proposition \ref{prop:cutsystem} that $N' \subseteq N$ induces isomorphisms in $H^2$ and $H^0$, so in particular we still have $H^2(N') = H^2(N) = 0$ and $H^0(N') = H^0(N) = \mathbb{Z}_2$. It follows easily from Lefschetz duality that $\partial N'$ has the homology of a $2$--sphere, and so it is a $2$--sphere. Since $\partial N'$ is polyhedral, it follows from the polyhedral Sch\"onflies theorem of Alexander (\cite[Theorem 12, p. 122]{moise2}) that $N'$ is a $3$--ball. Thus $N$ is a handlebody and $\{D_i\}$ is a complete system of cutting disks for $N$ satisfying the required conditions.
\end{proof}

\subsection{} \label{subsec:coh_generate} We now prove Lemma \ref{lem:kernel}. This requires some preliminary work.
\medskip

(i) All the inclusions \[(N,K) \subseteq (N,N^+) \supseteq (N^i \cup N^+, N^+) \supseteq (N^i,n^+)\] induce isomorphisms in \v{C}ech cohomology. The last one is just the strong excision property of \v{C}ech cohomology, while the other two make use of the flow to construct appropriate ``infinite time'' deformation retractions. For example, consider the map $H : (N,N^+) \times [0,+\infty) \longrightarrow (N,N^+)$ given by $(p,t) \longmapsto p \cdot \max \{-t,t^i(p)\}$. Fixing $p$ and allowing $t$ to go from $0$ to $+\infty$, the point $H(p,t)$ just follows the trajectory of $p$ backwards until it first hits $N^i$ (if it ever does) and remains stationary thereafter. The sets $H(N \times [k,+\infty))$ for $k = 0,1,2,\ldots$ form a nested sequence of compacta and the inclusion of each of them in the previous one is a homotopy equivalence, with $H|_{N \times [k,k+1]}$ providing a homotopy inverse. The first member of the sequence ($k = 0$) is $(N,N^+)$ and the limit of the sequence (the intersection of the whole sequence) is $(N^i \cup N^+,N^+)$. Thus the inclusion $(N^i \cup N^+, N^+) \subseteq (N,N^+)$ induces isomorphisms in \v{C}ech cohomology by its continuity property (\cite[Theorem 3.1, p. 261]{eilenbergsteenrod}). A similar argument applies to the inclusion $K \subseteq N^+$, and then the five lemma ensures that $(N,K) \subseteq (N,N^+)$ also induces isomorphisms in \v{C}ech cohomology.

(ii) Let $P \subseteq N^i$ be a compact $2$--manifold which is a neighbourhood of $n^+$. Then $H^k(N^i,P) \cong H^k(\overline{N^i \setminus P},\partial P) \cong H_{2-k}(\overline{N^i \setminus P},\partial N^i)$ where we have used excision and Lefschetz duality (in the manifold $\overline{N^i \setminus P}$), respectively. The geometric interpretation of this isomorphism (in dimension $1$, for example) is as follows. Let $z_{\gamma} \in H^1(N^i,P)$ correspond to $\gamma \in H_1(\overline{N^i \setminus P},\partial N^i)$. Then $z_{\gamma}$ acts on any $\alpha \in H_1(N^i,P)$ by returning its intersection number with $\gamma$ (modulo $2$). That is, if $c$ and $a$ are transverse polygonal $1$--chains representing $\gamma$ and $\alpha$, then $z_{\gamma}(\alpha)$ is the number of intersection points of $c$ and $a$ modulo $2$. This description characterizes $z_{\gamma}$ completely (by the universal coefficients theorem). Passing to the limit as $P$ gets smaller and using the continuity property of \v{C}ech cohomology and the compact supports property of homology yields an isomorphism $H^k(N^i,n^+) \cong H_{2-k}(N^i - n^+,\partial N^i)$. The geometric interpretation of this isomorphism is essentially the same as above.
\medskip

\begin{lemma} Assume $H^2(N,K) = 0$. Let $\gamma$ run over all arcs properly embedded in $N^i - n^+$. Then the $z_{\gamma}$ generate $H^1(N^i,n^+)$.
\end{lemma}
\begin{proof}  Consider the homology group $H_1(N^i - n^+,\partial N^i)$. Assuming everything to be triangulated, any homology class has a representative $c$ which is a (finite) sum of simplices whose algebraic boundary $\partial c$ lies in $\partial N^i$. Since we are taking coefficients in $\mathbb{Z}_2$ the chain $c$ can be identified with its underlying set, and it is then very easy to prove (using $\partial c = 0$) that the simplices in $c$ can be grouped into a structure $\sum c_i + \sum a_j$ where each $c_i$ is a simple closed curve and each $a_j$ is an arc whose endpoints lie in $\partial N^i$.

Since $H^2(N,K) = 0$, from (ii) above we have $H_0(N^i-n^+,\partial N^i) = H^2(N,K) = 0$ and so each component of $N^i - n^+$ intersects $\partial N^i$. Then each summand $c_i$ is homologous to a properly embedded arc; it suffices to join $c_i$ to $\partial N^i$ with a thin ribbon intersecting $c_i$ and $N^i$ along its short sides and add the boundary of the ribbon to $c_i$ (perhaps one needs to push $c_i$ slightly into the interior of $N^i$ first). Thus we have shown that $H_1(N^i - n^+,\partial N^i)$ is generated by (the homology classes defined by) properly embedded arcs $\gamma$, and so by (ii) above $H^1(N^i,n^+)$ is generated by the $z_{\gamma}$.
\end{proof}

Recall from (i) that the inclusion induced map $j^* : H^1(N,N^+) \longrightarrow (N^i,n^+)$ is an isomorphism. It follows from the previous lemma that $\{(j^*)^{-1}(z_{\gamma})\}$ is a generating system for $H^1(N,N^+)$. These cohomology classes have a simple geometric interpretation. We discussed in Subsection \ref{subsec:prep} how every arc $\gamma$ gives rise via the flow to a disk $D$ in $N$. There is a cohomology class $z_D$ that counts intersections (modulo $2$) with the disk $D$; evidently the restriction of this to $N^i$ counts intersections with $D \cap N^i = \gamma$; in other words, it is $z_{\gamma}$. Thus $(j^*)^{-1}(z_{\gamma}) = z_D$. Finally, since the inclusion $(N,K) \subseteq (N,N^+)$ also induces isomorphisms in \v{C}ech cohomology, we have:

\begin{lemma} Assume $H^2(N,K) = 0$. Let $\gamma$ run over all arcs properly embedded in $N^i - n^+$ and let $D$ run over all the disks generated from these arcs by the flow. Then $\{z_D\}$ generates $H^1(N,K)$.
\end{lemma}

Denote by $w_D$ the image of $z_D$ under the inclusion induced homomorphism $H^1(N,K) \longrightarrow H^1(N)$. Again, $w_D$ just counts intersections with $D$. By the previous lemma and the exactness of the cohomology sequence for the pair $(N,K)$, the set $\{w_D\}$ generates the kernel of the inclusion induced homomorphism $H^1(N) \longrightarrow H^1(K)$. This proves Lemma \ref{lem:kernel}.

\section{Proof of Theorem \ref{thm:main2} \label{sec:proof2}}

Very roughly, the proof of Theorem \ref{thm:main2} from the Introduction can be described as follows. First one shows that for any coloured $3$--ball $B$ there is a vectorfield on (a small neighbourhood of) the ball which generates a flow that realizes $B$ as an isolating block whose maximal invariant subset is a single rest point. Then to prove Theorem \ref{thm:main2} one starts with the coloured handlebody $N$, cuts it along the disks $\{D_i\}$ to obtain a $3$--ball $B = N \setminus \bigcup D_i \times (-1,1)$, finds a vectorfield as just described on $B$, and then extends it to the $1$--handles $D_i \times [-1,1]$ in such a way that no new invariant structure is introduced.

The idea is simple enough and the constructions are not difficult to visualize, but the formalization is slightly cumbersome. We will discuss a simpler case first (Proposition \ref{prop:realize}) to introduce the basic elements of the construction and then elaborate on it in Proposition \ref{prop:realize_handles}. In these two results $N$ will be any tame coloured manifold in $\mathbb{R}^3$, not necessarily a handlebody, since the proof is the same.

We would like to use the convenient language of vectorfields to define the flow. In order to do this it will be convenient to use the following trick. Recall that any $3$--manifold has a differentiable structure. We are given the coloured manifold $N$ as a subset of $\mathbb{R}^3$, but by regarding it as an abstract $3$--manifold we may assume that it is differentiable. We will enlarge it by attaching a collar $\partial N \times [0,1]$ on its boundary to obtain $M = N \cup \partial N \times [0,1]$ and define a flow $\varphi$ on $M$ which realizes $N$ as an isolating block. The advantage is that on $M$ we can use the language of differential topology to define the flow relatively easily. Returning from this abstract flow to $\mathbb{R}^3$ is simple. Denote by $N_{\rm emb}$ (for ``embedded'') the $3$--manifold $N$ thought of as a subset of $\mathbb{R}^3$. It is here that the assumption that $N_{\rm emb}$ be tame will come into play: it ensures that $N_{\rm emb}$ has a collar in $\mathbb{R}^3$, and so it can be slightly enlarged to $M_{\rm emb} = N_{\rm emb} \cup \partial N_{\rm emb} \times [0,1]$. The (identity) homeomorphism $N \longrightarrow N_{\rm emb}$ can be trivially extended to a homeomorphism $M \longrightarrow M_{\rm emb}$ and then the flow $\varphi$ defined on $M$ can be copied to a flow defined on $M_{\rm emb}$ which realizes $N_{\rm emb}$ as an isolating block. The flow can then be extended to all of $\mathbb{R}^3$ in whatever way.

\subsection{}

Let $N$ be a compact $3$--manifold. A spine of $N$ is a compact set $K \subseteq {\rm Int}\ N$ such that $N \setminus K$ is homeomorphic to $\partial N \times (-1,0]$, where $\partial N \ni x \longmapsto (x,0)$ under the homeomorphism. For example, a point in the interior of a $3$--ball is a spine.

 \begin{proposition} \label{prop:realize} Let $(N,P,Q)$ be a tame coloured manifold in $\mathbb{R}^3$ and $K$ a spine for $N$. There exists a flow $\varphi$ in $\mathbb{R}^3$ which realizes $N$ as an isolating block whose maximal invariant subset is $K$.
 \end{proposition}
 \begin{proof} Figure \ref{fig:def_V0}.(a) shows the annulus $A := \{u \in \mathbb{R}^2 : 1 \leq \|u\| \leq 3\}$ coloured in a certain way. On $A$ we consider the tangent vectorfield $u \longmapsto (\|u\|-1)(3-\|u\|)u$. The flow it generates has both boundary components comprised of fixed points and otherwise its trajectories evolve radially from the white component of the boundary towards the gray one.
\medskip

\fbox{Step 1} Consider $\partial N$ as a surface in the abstract. We may assume that it is differentiable and so are its subsets $P$ and $Q$. 
 We define a tangent vectorfield $W_0$ on $\partial N$. For each tangency curve $\tau$ let $A_{\tau}$ be a thin annulus along $\tau$. Each $A_{\tau}$ is diffeomorphic to the model $A$ of Figure \ref{fig:def_V0}.(a) via a colour preserving diffeomorphism. We define $W_0$ on each annulus $A_{\tau}$ by using any such diffeomorphism to copy the radial vectorfield on $A$. We will also use the expressions ``radial segment'' or ``radial direction'' on $A_{\tau}$ with the obvious meaning.
 
 Deleting the interiors of the annuli $A_{\tau}$ from $P$ and $Q$ we obtain slightly shrinked copies of $P$ and $Q$ that we shall denote by $sP$ and $sQ$. Define $W_0$ to be zero on $sP$ and $sQ$. The flow on $\partial N$ generated by $W_0$ has every point in $sP$ and $sQ$ stationary and otherwise flows radially across the annuli $A_{\tau}$, from $sP$ towards $sQ$.

Finally, define a smooth map $\theta : \partial N \longrightarrow [-1,1]$ such that
\begin{itemize}
	\item[(i)] $\theta|_{sP} \equiv +1$, $\theta|_{sQ} \equiv -1$ and $\theta|_{\tau} \equiv 0$ for each $t$-curve $\tau$;
	\item[(ii)] $\theta$ is strictly decreasing along each of the radial segments that fiber the annuli $A_{\tau}$.
\end{itemize}
For example, $\theta|_{A_{\tau}}$ could be (a suitably rescaled version of) the radius function of the annulus $A_{\tau}$.

\begin{figure}[h]
    \begin{pspicture}(0,0)(7,7)
	%\psgrid(0,0)(7,7)
	\rput[bl](0,0){\scalebox{0.5}{\includegraphics{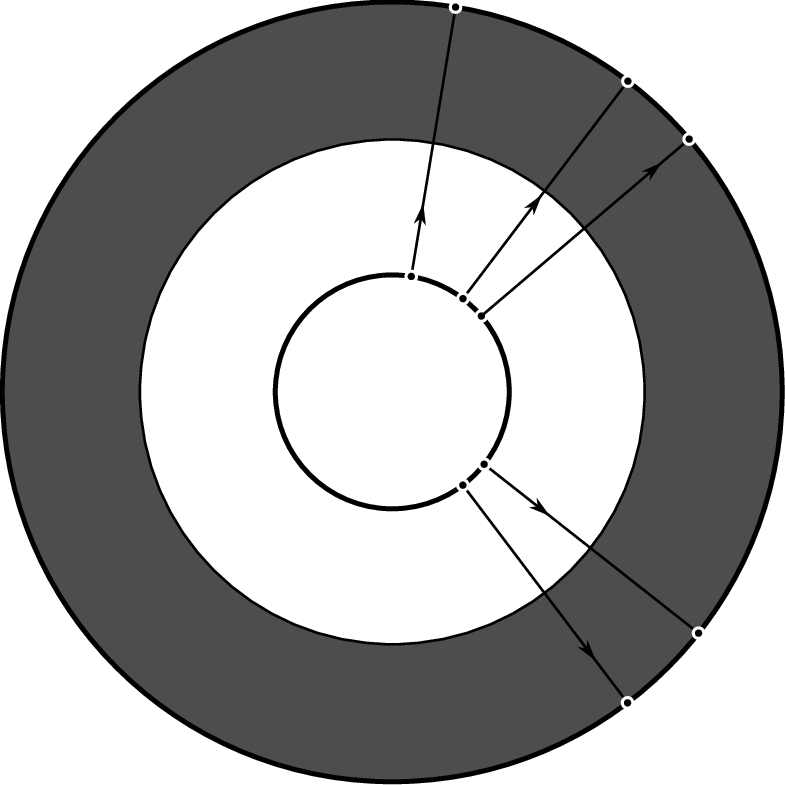}}}
        \rput[tr](0.9,0.9){$A$}
        %\rput[bl](1.9,1.9){$\tau$}
    \end{pspicture}
\caption{\label{fig:def_V0}}
\end{figure}

\medskip

\fbox{Step 2} Consider the abstract Cartesian product $\partial N \times [-1,1]$. We denote its points by $(x,r)$ and identify its tangent space at a point $(x,r)$ with $T_x (\partial N) \oplus \mathbb{R}$. Thus a tangent vector $V = (W,U)$ at a point $(x,r)$ consists of two components: a vector $W$ tangent to $\partial N$ at $x$ and a real number $U$.

Define a tangent vectorfield $V = (W,U)$ on $\partial N \times [-1,1]$ by \[W_{(x,r)} =  W_0(x) \ \ \text{ and } \ \ U_{(x,r)} = (r^2-1)\, \theta(x).\] For $r = \pm 1$ we have $V = (W_0,0)$ and so $V$ is actually tangent to the boundary $\partial N \times \{\pm 1\}$ of $\partial N \times [-1,1]$. Thus we can integrate $V$ to obtain a complete flow $\psi$ on $\partial N \times [-1,1]$.

If we take a cross section of $\partial N \times [-1,1]$ along one of the radial fibers of an annulus $A_{\tau}$ and extending a little into $sP$ and $sQ$ we observe the phase portrait shown in Figure \ref{fig:fases}. The gray half of the picture is $Q \times [-1,1]$ and the white half is $P \times [-1,1]$. Horizontal and vertical movement are controlled by the components $W$ and $U$ of $V$ respectively. 

\begin{figure}[h]
\begin{pspicture}(0,0)(15,6)
	%\psgrid(0,0)(15,6)
	\rput[bl](1.1,1){\scalebox{0.7}{\includegraphics{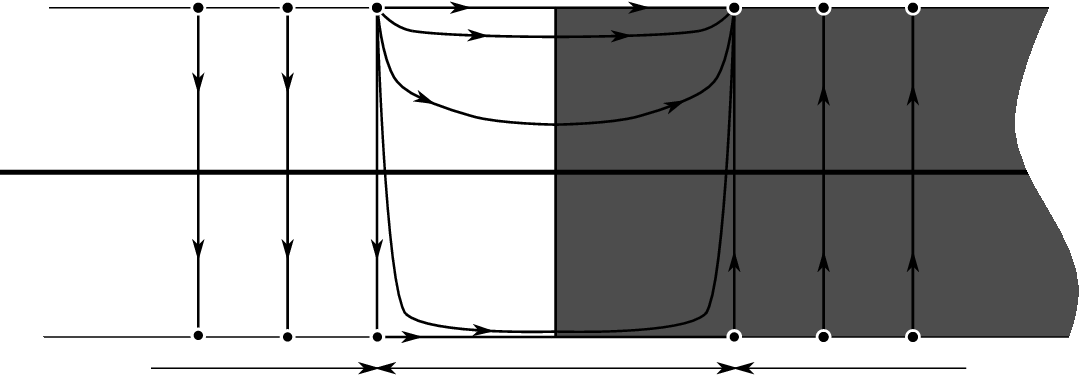}}}
	\rput[r](0.8,1.45){$r = -1$}
	\rput[r](0.8,3.4){$r =  0$}
	\rput[rb](0.8,5.25){$r =  1$}
        \rput[t](7.7,0.8){$A_{\tau}$}
        \rput[t](4.4,0.8){$sP$} \rput[t](11,0.8){$sQ$}
\end{pspicture}
\caption{\label{fig:fases}}
\end{figure}

We justify briefly the drawing of the phase portrait. On the top and bottom lines $r = \pm 1$ we have $V = ( W_0,0)$ so $\psi$ behaves there as explained in Step 1. Trajectories in $sP \times [-1,1]$ and $sQ \times [-1,1]$ move vertically because $\theta = 0$ there. Trajectories in $({\rm Int}\ A_{\tau}) \times (-1,1)$ have the U-shape suggested by the diagram for the following reason. Let $t \longmapsto (x(t),r(t))$ be one such trajectory and observe that $x(t)$ is an integral curve (in $\partial N$) for $W_0$ because $W$ only depends on $x$. Thus $x(t)$ must travel from the white boundary component of an annulus $A_{\tau}$ towards the other. While it does so $\theta(x(t))$ is positive until $x(t)$ hits $\tau$, where $\theta$ vanishes, and becomes negative thereafter. Since $\dot{r}(t)$ is proportional to $\theta(x(t))$ through a negative factor, it follows that $r(t)$ strictly decreases, reaching a minimum as $x(t)$ crosses $\tau$, and then strictly increases again. This justifies the U-shape of the trajectory. Notice also that for $t \rightarrow \pm \infty$ we must have $r(t) \rightarrow 1$.

Consider $\partial N \times (-1,0] \subseteq \partial N \times [-1,1]$. It follows from the previous paragraph that its maximal invariant subset is empty. Notice that along its boundary $\partial N \times \{0\}$ the flow $\psi$ behaves in the manner required for an isolating block: at points $(x,0)$ with $x \in {\rm Int}\ P$ we have $U_{(x,0)} = -\theta(x) < 0$ so that $\psi$ crosses $\partial N \times \{0\}$ transversally into $\partial N \times (-1,0]$, with the opposite happening when $p \in {\rm Int}\ Q$. For points $x$ in $\bigcup \tau \times \{0\}$ the flow $\psi$ is externally tangent to $\partial N \times \{0\}$.
\medskip

\fbox{Step 3} Still considering $N$ as an abstract $3$--manifold, enlarge it to a manifold $M$ by attaching a collar $\partial N \times [0,1]$ to its boundary via the map $(x,0) \longmapsto x$. Since $K$ is a spine of $N$, there exists a homeomorphism $\partial N \times (-1,0] \cong N \setminus K$; this can be extended to the collar in an obvious manner yielding a homeomorphism $\partial N \times (-1,1] \cong M \setminus K$. Actually, since any $3$--manifold has a unique differentiable structure, there is a diffeomorphism $\partial N \times (-1,1] \cong M \setminus K$. We may thus identify these two manifolds and think of $V$ as a vectorfield in $M \setminus K$. Define a vectorfield $\tilde{V}$ on all of $M$ by picking a continuous function $\mu : M \longrightarrow [0,1]$ such that $K = \mu^{-1}(0)$ and setting $\tilde{V} = 0$ over $K$ and $\tilde{V} = \mu V$ outside $K$. This vectorfield generates a flow $\varphi$ on $M$ which is stationary on $K$ and has the same phase portrait of Figure \ref{fig:fases} on $M \setminus K$. Thus $N$ is an isolating block in the phase space $M$, with entry and exit sets $P$ and $Q$ respectively, and maximal invariant subset $K$. We finally return to the embedded $N \subseteq \mathbb{R}^3$ as explained at the beginning of this section.
\end{proof}

\subsection{}

 Suppose $(N,P,Q)$ is a (tame) coloured manifold in $\mathbb{R}^3$ and $\{D_i\}$ is a system of cutting disks for $N$. Let $N'$ be the result of cutting $N$ along the $\{D_i\}$ and let $K' \subseteq {\rm Int}\ N'$ be a spine of $N'$.

 \begin{proposition} \label{prop:realize_handles} There exists a flow $\varphi$ in $\mathbb{R}^3$ which realizes $N$ as an isolating block whose maximal invariant subset is $K'$.
 \end{proposition}

Notice the difference with Proposition \ref{prop:realize}: here the maximal invariant subset of $N$ is not a spine $K$ of $N$ itself, but a spine $K'$ of the simpler manifold $N'$.

\begin{proof}[Proof of Proposition \ref{prop:realize_handles}] For notational simplicity we prove the proposition when there is just one cutting disk $D$. The proof is a modification of the argument of Proposition \ref{prop:realize}.
\medskip

Let $N'$ be the result of cutting $N$ along the disk $D$. As described before, $\partial N'$ is coloured almost completely save for the two lined cutting disks $D_{\pm}$. For the purpose of comparison with Proposition \ref{prop:realize}, extend the colouring of $\partial N'$ to the disks $D_{\pm}$ in the obvious manner shown in Figure \ref{fig:1-handle_col}, painting one half of each disk black and the other white. Thus now $N'$ is a coloured manifold whose $t$--curves we denote generically by $\tau'$. Notice that each of the disjoint disks $D_{\pm}$ intersect $\bigcup \tau'$ along a diameter. For later reference we label by $t_-$ and $t_+$ the points shown in Figure \ref{fig:1-handle_col}. The point $t_-$ is any one of the two points in $\partial D_-$ where the colouring changes; once this is chosen, $t_+$ is the point in $\partial D_+$ which sits ``in front'' of $t_-$; i.e. which is joined to $t_-$ by a portion of a $t$--curve in $D \times [-1,1]$.

\begin{figure}[h!]
\begin{pspicture}(0,0)(5,5)
	%\psgrid(0,0)(5,5)
	\rput[bl](0,0){\scalebox{0.8}{\includegraphics{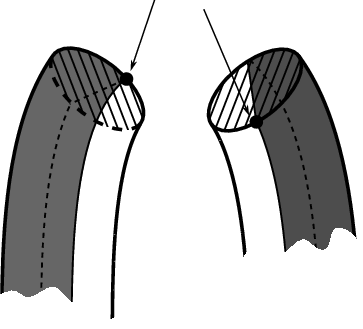}}}
        \rput[t](4.8,1){$N'$}
        \rput[b](0.8,3.8){$D_-$}
        \rput[bl](4,3.8){$D_+$}
        \rput[b](2.1,4.5){$t_-$} \rput[b](2.8,4.4){$t_+$}
\end{pspicture}
\caption{\label{fig:1-handle_col}}
\end{figure}

\medskip

\fbox{Step 1} We first define a tangent vectorfield $W_0$ on $\partial N'$. Again, for each tangency curve $\tau'$ of $\partial N'$ we let $A_{\tau'}$ be a thin annulus along $\tau'$. We use the same construction as in Step 1 of Proposition \ref{prop:realize} with the following modification. In the boundary of $N'$, somewhere along $\bigcup \tau'$, there are the two cutting disks $D_{\pm}$. We choose the annuli $A_{\tau'}$ in such a way that the disks $D_{\pm}$ are contained in these annuli. Also, when choosing the diffeomorphisms $A_{\tau'} \cong A$ to copy the radial vectorfield in $A$ to a vectorfield $W_0$ in the annuli $A_{\tau'}$, we take the following precaution. Suppose that (say) $D_+$ is contained in $A_{\tau'}$. We then choose the diffeomorphism $A_{\tau'} \cong A$ in such a way that $D_+ \subseteq A_{\tau'}$ corresponds to the lined disk $D \subseteq A$ in Figure \ref{fig:def_V0_disc}.(a). We define $D'_+$ to be the preimage of $D'$ under the diffeomorphism $A_{\tau'} \cong A$. Finally, we extend $W_0$ by zero outside the annuli $A_{\tau'}$.

For later reference we observe the following. Consider the phase portrait of Figure \ref{fig:def_V0_disc}.(b) and denote by $\mathbb{D}_1 \subseteq \mathbb{D}_2$ the two concentric disks depicted there. Every point in $\partial \mathbb{D}_2$ is a rest point and otherwise trajectories flow in a parallel fashion. The point $t \in \partial \mathbb{D}_1$ is one of the two points where a trajectory is tangent to $\mathbb{D}_1$. Now imagine one modifies the radial flow in $A$ in panel (a) of the figure by stopping it outside the interior of $D'$. Every point in the boundary of $D'$ becomes a rest point and trajectories inside $D'$ still flow radially as in the drawing. This modified flow on $D'$ is conjugate to the model shown in Figure \ref{fig:def_V0_disc}.(b) via a colour preserving homeomorphism $(D',D) \cong (\mathbb{D}_2,\mathbb{D}_1)$. In particular the two points in $\partial D$ where the modified radial flow in $D'$ is tangent to $D$ must go to the corresponding tangency points in $\partial \mathbb{D}_1$, and one can choose which one goes to $t$.

\begin{figure}[h]
\null\hfill
\subfigure[]{
    \begin{pspicture}(0,0)(7,7)
	%\psgrid(0,0)(7,7)
	\rput[bl](0,0){\scalebox{0.5}{\includegraphics{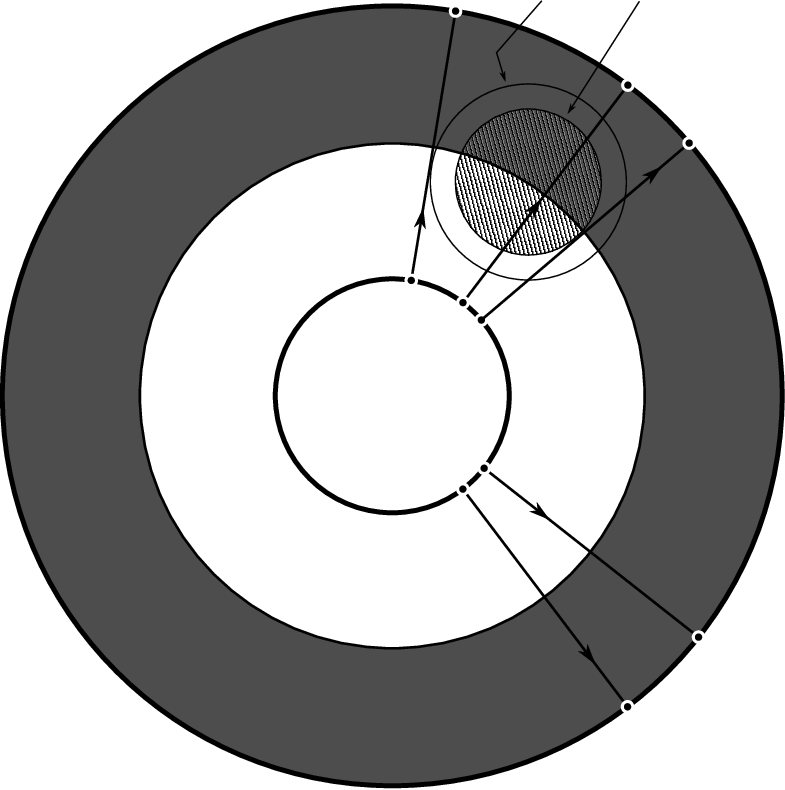}}}
        \rput[tr](0.9,0.9){$A$}
        \rput[bl](4.6,6.7){$D'$}
        \rput[bl](5.4,6.7){$D$}
    \end{pspicture}}
\hfill
\subfigure[]{
    \begin{pspicture}(-0.2,0)(4,7)
	%\psgrid(-0.2,0)(4,7)
	\rput[bl](0,1.5){\scalebox{0.5}{\includegraphics{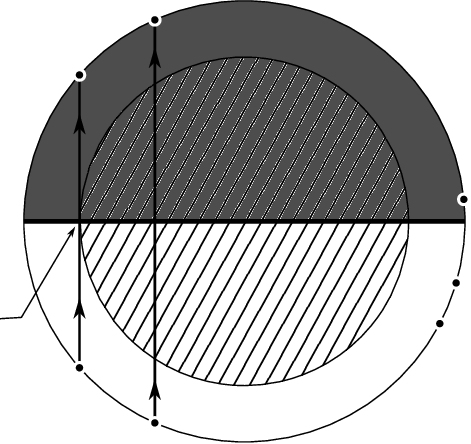}}}
    \rput[r](-0.1,2.6){$t$}
    \rput[tl](3.1,1.7){$\mathbb{D}_2$}
    \rput(2.8,2.8){$\mathbb{D}_1$}
    \end{pspicture}}
\hfill\null
\caption{\label{fig:def_V0_disc}}
\end{figure}

The map $\theta$ is defined exactly as before.
\medskip

\fbox{Step 2} Now we define a tangent vectorfield $V = (W,U)$ on $\partial N' \times [-1,1]$ using $W_0$ in a manner similar to the previous proposition, but with some modifications so that we can later on extend the dynamics to all of $N$. Let $\alpha : \partial N' \longrightarrow [0,1]$ be a smooth map which is strictly positive exactly on the interior of the disks $D'_{\pm}$. Set \[W_{(x,r)} = ((1+r) \alpha(x) + 1 - r) \cdot W_0(x) \ \ \text{ and } \ \ U_{(x,r)} = (r^2-1)\, \theta(x).\] As before, $V$ generates a complete flow $\psi$ on $\partial N' \times [-1,1]$. We observe the following:

\begin{itemize}
    \item[(i)] For $r = 1$ we have $V = (2 \alpha W_0,0)$. Every point outside the interior of the disks $D'_{\pm} \times \{1\}$ is an equilibrium for $\psi'$, and on the interior of the disks themselves the orbits evolve in the usual radial fashion. Thus as discussed before $\psi|_{D'_{\pm} \times \{1\}}$ is conjugate, via a colour preserving homeomorphism, to the model of Figure \ref{fig:def_V0_disc}.(b).
    \item[(ii)] In the region $-1 \leq r < 1$ the phase portrait of $\psi$ is qualitatively the same as that of Figure \ref{fig:fases}. The reason is the following. Suppose $t \longmapsto (x(t),r(t))$ is an integral curve of $V$ passing through a point $(x_0,r_0)$ with $x_0$ in the interior of some annulus $A_{\tau'}$ and $r_0 < 1$. Since the subset $r = 1$ is invariant, we must have $r(t) < 1$ for all $t$, and in particular the numerical factor which multiplies $W_0(x)$ in the definition of $W_{(x,r)}$ is strictly positive for all $t$. Now, the derivative of  $t \longmapsto \theta(x(t))$ is given by said numerical factor times the derivative of $\theta$ in the direction $W_0(x(t))$. The latter is strictly negative by definition of $\theta$, and so $\theta(x(t))$ is strictly decreasing. Thus again $x(t)$ evolves monotonically across the annulus $A_{\tau'}$, from the white boundary towards the gray boundary. Accordingly $r(t)$ strictly decreases until $\theta = 0$ and then strictly increases again. 
\end{itemize}

Figure \ref{fig:fases_asa}.(a) shows again a cross section of the flow along a radial segment of an annulus $A_{\tau'}$ that intersects the disk $D'_+$ (say). Notice how points with $r = 1$ outside $D'_+$ are now equilibria while inside $D'_+$ we still have the same radial flow as before.

\begin{figure}
\null\hfill
\subfigure[]{
\begin{pspicture}(0,0)(6,5)
    %\psgrid(0,0)(6,5)
    \rput[bl](0,0){\scalebox{0.7}{\includegraphics{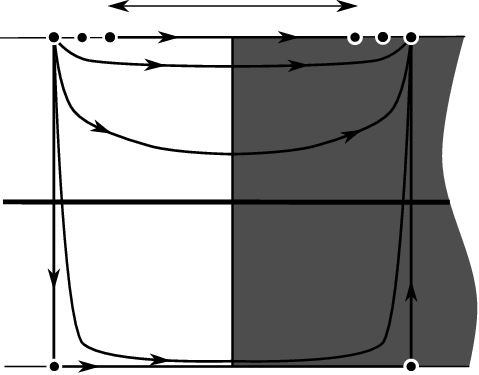}}}
    \rput[b](2.8,4.4){$D'_+$}
\end{pspicture}}
\hfill
\subfigure[]{
\begin{pspicture}(0,0)(6,5)
    %\psgrid(0,0)(6,5)
    \rput[bl](0,0){\scalebox{0.7}{\includegraphics{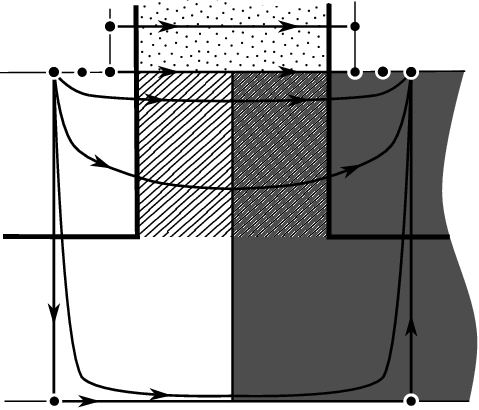}}}
\end{pspicture}}
\hfill\null
\caption{\label{fig:fases_asa}}
\end{figure}
\medskip

\fbox{Step 3} As in Step 3 of Proposition \ref{prop:realize}, enlarge $N'$ to a manifold $M'$ by attaching a collar onto its boundary, identify $\partial N' \times (-1,1] \cong M' \setminus K'$ and think of $V$ as a vectorfield on $M' \setminus K'$. Multiply it by an appropriate function $\mu$ to extend it continuously by zero to $K'$ and use this to generate a new flow $\varphi_{M'}$ in $M'$. This flow has $N'$ as an isolating neighbourhood for $K'$ and its entry and exit pattern into $N'$ accords to the colouring of $\partial N'$ described at the beginning of the proof. Moreover, $\varphi_{M'}$ sweeps across the disks $D'_{\pm} \times \{1\} \subseteq \partial M'$ as in Figure \ref{fig:def_V0_disc}.(b). The goal now is to attach back to $N'$ the $1$--handle we removed when cutting along the disk $D$, while at the same time extending the flow $\varphi_{M'}$ appropriately. This we do as follows.

Consider again the phase portrait of Figure \ref{fig:def_V0_disc}.(b). Define $H := \mathbb{D}_2 \times [-1,1]$ and endow it with the flow $\varphi_H$ which repeats Figure \ref{fig:def_V0_disc}.(b) on each slice $\mathbb{D}_2 \times \{t\}$. One should think of $H$ as a $1$--handle with pasting disks $\mathbb{D}_2 \times \{\pm 1\}$ which we are going to attach onto $M'$ along the disks $D'_{\pm} \times \{1\} \subseteq \partial N' \times \{1\} \subseteq \partial M'$. This we do by choosing pasting homeomorphisms $\mathbb{D}_2 \times \{\pm 1\} \cong D'_{\pm} \times \{1\}$ which: (a) conjugate $\varphi_H$ and $\varphi_{M'}$, (b) are colour preserving, (c) send $\mathbb{D}_1 \times \{\pm 1\}$ onto $D_{\pm} \times \{1\}$, (d) send $(t_{\pm},1)$ to $(t,\pm 1)$. These homeomorphisms exist by the discussion at the end of Step 1. The resulting space $M := M' \cup H$ is a $3$--manifold (with corners, but we do not need to worry about differentiability any more), and the flows on $M'$ and on $H$ match up on the pasting disks so there is a well defined flow $\varphi$ on $M$.

Consider the solid tube $\mathbb{D}_1 \times [-1,1] \subseteq H$ and lengthen it slightly to the solid tube \[T := D_+ \times [0,1] \cup \mathbb{D}_1 \times [-1,1] \cup D_- \times [0,1] \subseteq M.\] The bases of this tube are the disks $D_{\pm} \times \{0\}$, which lie on $\partial N'$. Figure \ref{fig:fases_asa}.(b) shows a cross sectional view of this. The lined region is $D_+ \times [0,1]$ whereas the dotted region extending upwards is (a portion of) $\mathbb{D}_1 \times [-1,1]$; the other end of the tube $T$ looks similar. It is clear from this that the whole tube $T$ is coloured by the flow $\psi$ as shown in Figure \ref{fig:1-handle}; one should think of the flow as running vertically from the bottom up in the figure.

\begin{figure}[h!]
\begin{pspicture}(-0.4,0)(8,4)
	%\psgrid(-0.4,0)(8,4)
	\rput[bl](0,1){\includegraphics{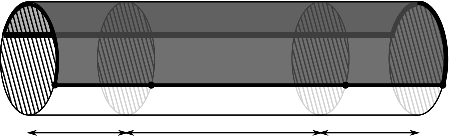}}
        \rput[b](1,3.5){$T$}
	\rput[t](0.5,0.8){$D_+ \times [0,1]$}
        \rput[t](7.2,0.8){$D_- \times [0,1]$}
        \rput[b](4,2.1){$a$}
        \rput[t](3.8,0.8){$\mathbb{D}_1 \times [-1,1]$}
        \rput[t](2.6,1.75){$t$} \rput[t](5.9,1.75){$t$}
        \rput[t](7.7,1.75){$t_-$} \rput[t](1.1,1.75){$t_+$}
\end{pspicture}
\caption{\label{fig:1-handle}}
\end{figure}

It should be clear from the construction so far that $N' \cup T \subseteq M := M' \cup H$ is an isolating block for $K'$ having $\mathbb{D}^1 \times \{0\}$ as a cutting disk. Notice that there are two arcs of tangency points running along the lateral face of $T$. Recall that at the beginning of the proof we distinguished two points $t_{\pm}$ in $\partial D_{\pm}$ as well as a point $t \in \partial \mathbb{D}_1$ and required that when pasting $H$ onto $M'$ the points $(t_+,1)$ and $(t,-1)$ matched, and similarly at the other end of $H$. This ensures that the two points $(t_{\pm},0) \in \partial N'$ cobound one of the tangency arcs along $T$; namely \[a = \underbrace{\{t_+\} \times [0,1]}_{\subseteq D_+ \times [0,1]} \cup \underbrace{\{t\} \times [-1,1]}_{\subseteq \mathbb{D}_1 \times [-1,1]} \cup \underbrace{\{t_-\} \times [0,1]}_{\subseteq D_- \times [0,1]}.\] See Figure \ref{fig:1-handle}.
\medskip

\fbox{Step 4} Finally we return back to the original problem in $\mathbb{R}^3$. As usual, we denote by $N_{\rm emb} \subseteq \mathbb{R}^3$ the original manifold $N$ embedded in $\mathbb{R}^3$. We have a colour preserving homeomorphism $h : N' \longrightarrow N'_{\rm emb}$ (essentially, the identity) and want to extend this to a colour preserving homeomorphism $h : N' \cup T \longrightarrow N_{\rm emb}$. Since $N_{\rm emb} = N'_{\rm emb} \cup D_{\rm emb} \times [-1,1]$, where $D_{\rm emb}$ is the disk along which we cut $N$ to begin with, we only need to extend $h$ to a colour preserving homeomorphism from $T$ to $D_{\rm emb} \times [-1,1]$. Both $T$ (Figure \ref{fig:1-handle}) and $D_{\rm emb} \times [-1,1]$ (Figure \ref{fig:cap}.(a)) are solid cylinders that look alike, including their colouring, so the existence of such an extension of $h$ seems plausible (there is a slight subtlety, though). To prove it we will use repeatedly the fact that a homeomorphism between the boundaries of two balls can be extended to a homeomorphism between the two balls (by coning from the center). 

Let us focus first on the two arcs of tangency points running along the boundary of $T$. One is $a$; denote the other by $b$. There are also two arcs of $t$--curve on the lateral face of $D_{\rm emb} \times [-1,1] \subseteq N_{\rm emb}$. One of them, say $a_{\rm emb}$, has $t_{\pm}$ as its endpoints by definition (back at the beginning of the proof). Denote the other by $b_{\rm emb}$. Now, $h$ is already defined on the lined disks $D_{\pm} \times \{0\} \subseteq \partial T$, and so in particular it is defined on the endpoints $(t_{\pm},0)$ of the arc $a$ which it carries precisely onto the endpoints of $a_{\rm emb}$. Thus it can be extended to carry the whole arc $a$ onto $a_{\rm emb}$, and similarly for $b$ and $b_{\rm emb}$. This is the subtlety we mentioned earlier: had we not been careful, $h$ might have interchanged the endpoints of $a$ and $b$, rendering a colour preserving extension impossible.

Now focus on the gray half of the lateral face of $\partial T$. It is a $2$--ball bounded by $a$, $b$, and half of each of $\partial D_{\pm} \times \{0\}$. The homeomorphism $h$ is already defined on all of them, and so we can extend it to a homeomorphism of the whole gray region onto the corresponding region of the lateral face of $D_{\rm emb} \times [-1,1]$. The same goes for the white half. We now have $h$ extended to the whole boundary of the $3$--ball $T$ onto the boundary of $D_{\rm emb} \times [-1,1]$, and so finally it can be extended to all of $T$. By construction this extended homeomorphism $h : N' \cup T \longrightarrow N_{\rm emb}$ is colour preserving. The proof then finishes in the usual way, by observing that $M = M' \cup H$ is $N$ with a collar attached and using this to extend $h$ to a homeomorphism of $M$ onto a collar of $N_{\rm emb}$ and copy the flow into $\mathbb{R}^3$.
\end{proof}

\subsection{} The proof of Theorem \ref{thm:main2} is now straightforward. We recall its statement:

\begin{theorem} Let $(N,P,Q) \subseteq \mathbb{R}^3$ be a tame coloured handlebody. Assume it has a complete system of cutting disks $\{D_i\}$ such that each $\partial D_i$ intersects $\bigcup \tau$ transversally at two points. Then there exists a flow $\varphi$ on $\mathbb{R}^3$ which realizes $N$ as an isolating block whose maximal invariant subset is a single point.
\end{theorem}
\begin{proof} The manifold $N'$ that results by cutting $N$ along the $\{D_i\}$ is a $3$--ball, which evidently has a spine $K'$ consisting of a single point. It follows directly from Proposition \ref{prop:realize_handles} that there is a flow $\varphi$ in $\mathbb{R}^3$ which realizes $N$ as an isolating block for the rest point $K'$.
\end{proof}

\section{Proof of Theorem \ref{thm:main3} \label{sec:proof3}}

The fundamental group of a handlebody $N$ is free in some generators $\{x_i\}$. When the handlebody is coloured, every $t$--curve defines an  element of this group (up to conjugacy since no basepoint is fixed); i.e. a word in the letters $x_i^{\pm 1}$ in the manner described in the Introduction. The geometric condition of Theorem \ref{thm:main1}.(ii) can be reformulated in terms of these words, and the powerful machinery of free groups applied to obtain Theorem \ref{thm:main3}. However, to keep the discussion as simple as possible we will hide these algebraic foundations of our discussion and adopt a very pedestrian approach instead, formulating everything in terms of purely syntactical manipulations of words.  

\subsection{} Let $\{x_1,x_1^{-1},\ldots,x_g,x_g^{-1}\}$ be a collection of letters, where each $x_i$ and $x_i^{-1}$ should be thought of as inverses of each other. A word $W$ means a (finite) sequence of letters. Its inverse $W^{-1}$ is the same sequence written backwards and with each letter replaced with its inverse. There is an empty word, denoted by $1$, with no letters. A word is cyclically reduced if it has no consecutive pairs of the form $x_i x_i^{-1}$ or $x_i^{-1} x_i$ and the same condition holds for its last and first letters (hence the word ``cyclically''). By successively cancelling pairs of letters of this sort until no further cancellation is possible, every word $W$ gives rise to a cyclically reduced word which is unique save for cyclic permutations (i.e. moving the last symbol of the word to the beginning of the word). The algebraic length (or simply ``length'') of a word is the number of letters in its cyclically reduced form.

Suppose $W$ is a word in the letters $x_i^{\pm 1}$. By substituting each appearance of $x_i$ in $W$ with some other word $V_i$, and those of $x_i^{-1}$ with $V_i^{-1}$, one performs a substitution. We consider substitutions $\tau$ of the following two types:
\begin{itemize}
    \item[(1)] $\tau$ permutes the letters $x_i^{\pm 1}$ among themselves.
    \item[(2)] For some fixed letter $a \in \{x_i^{\pm 1}\}$, called the multiplier, $\tau$ leaves $a^{\pm 1}$ unchanged and replaces every other $x_i$ with one of these four possibilities: $x_i$ itself, $x_i a$, $a^{-1} x_i$, or $a^{-1} x_i a$. In other words, each $x_i$ may or may not acquire an $a$ on the right and may or may not acquire an $a^{-1}$ on the left.
\end{itemize}

These are called Whitehead substitutions. Clearly they are all reversible, with the inverse being another Whitehead substitution. Also, there are only finitely many of these substitutions. For the sake of brevity let us write $S \sim S'$ to denote that $S'$ arises from $S$ by performing finitely many Whitehead substitutions sequentially, possibly accompanied by (cyclic) reductions. $\sim$ is an equivalence relation.

The algebraic length of a finite set of words is the sum of the algebraic lengths of its members. Fix some finite set of words $S$. We shall say that a set $S_{\rm min} \sim S$ is a minimal form of $S$ if it achieves the minimum algebraic length among all $S' \sim S$. There may exist several minimal forms of $S$, of course all with the same algebraic length. A celebrated theorem of Whitehead (\cite{whitehead1}; for a more condensed proof see \cite[Proposition 4.20, p. 35]{lyndonschupp}) states the following:

\begin{itemize}
    \item[(i)] If $S' \sim S$ is not minimal, there exists a Whitehead substitution which reduces its algebraic length.
    \item[(ii)] Given any two minimal forms of $S$, there exists a sequence of Whitehead substitutions that carries one onto the other while keeping the algebraic length the same at each step.
\end{itemize}

Notice in particular that (i) gives a procedure to find a minimal forms of $S$. One starts with $S$ and applies a Whitehead substitution to all the words in $S$, cyclically reducing the resulting words and obtaining a new set $S'$. If the algebraic length of $S'$ is not smaller than that of $S$, one tries again with a different Whitehead substitution. Since there are only finitely many of these, after some time one either finds a substitution $\tau$ which reduces the length of $S$ or finds that none does. In the latter case $S$ is already minimal. In the former case one throws $S$ away and starts again with $\tau(S)$. Evidently this process finishes after a finite number of steps (no more than the algebraic length of $S$) and returns a minimal set $S_{\rm min} \sim S$. We will call this algorithm ``Whitehead reduction''. See Example \ref{ex:whitehead_reduction} below for an illustration.

\subsection{} Given a set of words $S$ in the letters $x_i^{\pm 1}$, we are interested in the following algebraic condition:
\smallskip

$(A)$ For every $i$, the letters $x_i^{\pm 1}$ either do not appear at all among the words in $S$ or both appear, exactly once each.
\smallskip

Observe that if $S$ satisfies condition $(A)$, then after cyclically reducing its words it still satisfies the condition (since any reduction involves cancelling a pair $x_i^{\pm 1}$). Now we examine the effect of a Whitehead substitution:

\begin{lemma} \label{lem:Qinvariant} Suppose $S$ satisfies $(A)$. Let $\tau$ be a Whitehead substitution such that the length of $\tau(S)$ is not bigger than the length of $S$. Then, after cyclic reduction, $\tau(S)$ also satisfies $(A)$.
\end{lemma}
\begin{proof} We can take $S$ to be cyclically reduced. Otherwise we replace it with its cyclic reduction, which still satisfies $(A)$.

If $\tau$ is just a permutation of the variables the result is obvious, so suppose $\tau$ is a substitution with multiplier $a$. In passing from $S$ to $\tau(S)$ no new letters $x_i^{\pm 1}$ different from $a^{\pm 1}$ are inserted, so for those letters condition $(A)$ is satisfied in $\tau(S)$ because it was already satisfied in $S$. Thus to prove the lemma we only need to show that (after cyclic reduction) $a$ and $a^{-1}$ do not appear in $\tau(S)$, or appear exactly once each.

When one acts with $\tau$ on a word $W$ one obtains a new word $\tau(W)$ which is in principle longer than $W$ since several instances of $a$ and $a^{-1}$ will generally have been inserted. This might lead to the appearance of portions $...a a^{-1}...$ and $...a^{-1} a...$ that can be cancelled out. It can be shown that if the original word $W$ was cyclically reduced, only these cancellations involving the multiplier $a$ are possible in $\tau(W)$ (\cite[Proof of Proposition 4.16, p. 31]{lyndonschupp}). Thus the change in length from $S$ to $\tau(S)$ equals the change in the number of appearances of $a$ and $a^{-1}$ in $S$ to $\tau(S)$ (after cyclic reduction). In turn, this equals \[2(\# \text{appearances of } a \text{ in } \tau(S) - \# \text{appearances of } a \text{ in } S)\] because $a$ and $a^{-1}$ appear both the same number of times in $S$ and also in $\tau(S)$. Indeed, this property is true in $S$ for every letter by condition $(A)$ and is then preserved upon applying $\tau$, since whenever a letter $x_i$ acquires an $a$ on its right (say), $x_i^{-1}$ acquires an $a^{-1}$ on its left, balancing the count of $a$ versus $a^{-1}$ in $\tau(S)$ again.

Since by assumption the length of $\tau(S)$ is not bigger than that of $S$, it follows from the formula above that the number of appearances of $a$ (hence also of $a^{-1}$) in $\tau(S)$ is not bigger than its number of appearances in $S$, which is at most $1$ by condition $(A)$. This finishes the proof. 
\end{proof}

\begin{lemma} \label{lem:Qinvariant_min} Let $S$ be a finite set of words. Then the following are equivalent:
\begin{itemize}
    \item [(a)] There exists $S' \sim S$ which satisfies $(A)$.
    \item[(b)] Every minimal form of $S$ satisfies $(A)$.
\end{itemize}
In particular either every minimal form of $S$ satisfies $(A)$ or none does.
\end{lemma}
\begin{proof} Evidently (b) $\Rightarrow$ (a), since minimal forms always exist. To prove the converse assume $S' \sim S$ satisfies $(A)$ and $S_{\rm min}$ is any minimal form of $S$. By parts (i) and (ii) of Whitehead's theorem we can go from $S'$ to $S_{\rm min}$ by first applying a sequence of Whitehead substitutions $\tau$ that decrease the length at each step (until we arrive at some minimal form of $S$) and then another sequence of Whitehead substitutions $\tau$ that keeps the length constant at each step and arrives at $S_{\rm min}$. We know from Lemma \ref{lem:Qinvariant} that all these $\tau$ preserves condition $(A)$. Thus $S_{\rm min}$ satisfies $(A)$, because $S'$ did.
\end{proof}

In fact it can be shown that a set $S$ satisfying $(A)$ is already very close to being minimal. 

\subsection{} We finally return to the problem at hand; namely, recognizing when a coloured handlebody has a complete cut system that satisfies condition (ii) in Theorem \ref{thm:main1}. First we explain the geometric interpretation of the preceding discussion.

Let $N$ be a handlebody with a cut system $\{D_1,\ldots,D_g\}$. Cut $N$ open along this system to obtain a $3$--ball. Fix as a basepoint $*$ the center of the ball and consider $2g$ straight segments joining it with the center of each of the disks $D_i^{\pm}$ on the boundary of the ball. If we now glue together each pair $D_i^+, D_i^-$ we recover the handlebody $N$ together with a collection of $g$ closed curves $b_1,\ldots,b_g$ based at $*$ such that each of them pierces exactly one $D_i$, and just once. Their homotopy classes $\{[b_1],\ldots,[b_g]\}$ form a basis for the free group $\pi_1(N,*)$.

Suppose $s$ is an oriented closed curve in $\partial N$ and assign to it a word $V(x_1,\ldots,x_g)$ as described in the Introduction. Choose any path $\alpha$ that joins the basepoint $*$ to any point in $s$ and consider the concatenation $\alpha * s * \alpha^{-1}$, which is a loop based at $*$. The expression of $[\alpha* s* \alpha^{-1}]$ in the basis $\{b_i\}$ is then conjugate to the element $V([b_1],\ldots,[b_g]) \in \pi_1(N,*)$. This word depends on the choice of $\alpha$; however, its conjugacy class does not and so as a cyclically reduced word $V$ depends only on $s$. This is how cyclic words arise in our context.

Whitehead substitutions also have a geometric interpretation. Given a complete cut system $\{D_1,\ldots,D_g\}$ for a handlebody there is a geometric manipulation that generates a new system $\{D_1,\ldots,D'_i,\ldots,D_g\}$ by replacing one of the disks $D_i$ with a new disk $D'_i$ constructed as a band sum of the old $D_i$ and some other $D_j$. If an oriented curve $s$ in $\partial N$ spelled the word $V$ with respect to the cut system $\{D_i\}$, the word $V'$ that it spells with respect to $\{D'_i\}$ is obtained by a substitution with the general structure $x_j \rightarrow x_j x_i$ (with minor variations such as $x_j x_i^{-1}$, $x_i x_j$, etc. depending on the orientations assigned to the disks and how the band sum is performed). Observe that this is a Whitehead substitution. Now, it is a theorem that any two complete cut systems of $N$ can be related to each other by a sequence of these manipulations, and so the following lemma holds: 

\begin{lemma} \label{lem:changebasis} Suppose $N$ is a coloured handlebody of genus $g$ and $\{D_i\}$ and $\{D'_i\}$ are two cut systems for $N$. Let $S$ and $S'$ be the sets of words obtained by reading the $t$--curves with respect to these two cut systems. Then $S \sim S'$.
\end{lemma}

The above argument illuminates the geometric reason for the lemma but omits many details. Here is an alternative, algebraic proof:

\begin{proof}[Proof of Lemma \ref{lem:changebasis}] The two cut systems give rise to bases $\{[b_i]\}$ and $\{[b'_i]\}$ of $\pi_1(N)$ as explained before. Each $[b_i]$ can be expressed as $U_i([b'_1],\ldots,[b'_g])$ where the word $U_i(x_1,\ldots,x_g)$ records the intersections of $b_i$ with the $\{D'_i\}$. Let $V(x_1,\ldots,x_g)$ be one of the words in $S$; i.e. the word read off from one of the $t$--curves $\tau$ in $N$ with respect to the cut system $\{D_1,\ldots,D_g\}$. Then $[\tau] = V([b_1],\ldots,[b_g])$ in $\pi_1(N)$ and \[[\tau] = V(U_1([b'_1],\ldots,[b'_g]),\ldots,U_g([b'_1],\ldots,[b'_g]))\] so the word $V(U_1(x_1,\ldots,x_g),\ldots,U_g(x_1,\ldots,x_g))$ expresses $[\tau]$ in the basis $\{b_i\}$. Since such an expression is unique, the set of words $S'$ and the set of words $\{V(U_1,\ldots,U_g) : V \in S\}$ must coincide after cyclic reduction.

Since both $\{[b_i]\}$ and $\{[b'_i]\}$ are basis for $\pi_1(N)$, there exists an automorphism which carries the first onto the second (by the definition of a free group there exist endomorphisms of $\pi_1(N)$ that carry $\{[b_i]\}$ onto $\{[b'_i]\}$ and viceversa; their composition in both orders is the identity again by the definition of a free group). Whitehead substitutions encompass the so-called Nielsen substitutions, and the latter already suffice to generate all automorphisms of a free group (\cite[Proposition 4.1, p. 23]{lyndonschupp}). Thus there exists a sequence of Whitehead substitutions that carries $\{U_1,\ldots,U_g\}$ to $\{x_1,\ldots,x_g\}$. If we apply this sequence of Whitehead substitutions on the set $\{V(U_1,\ldots,U_g) : V \in S\}$ we evidently get (after cyclic reduction) the set $S$ back. Thus \[S' = \{V(U_1,\ldots,U_g): V \in S\} \sim S\] as claimed. 
\end{proof}

We can now prove Theorem \ref{thm:main2}, which we recall here as follows:

\begin{theorem} \label{thm:main3_restated} Let $N$ be a coloured handlebody with a nonempty collection of $t$--curves. Let $\{D_1,\ldots,D_g\}$ be any cut system for $N$ and let $S$ be the collection of words in the letters $x_i^{\pm 1}$ obtained by reading the $t$--curves, oriented as the boundary of the gray region. Denote by $S_{\rm min}$ any minimal form of $S$. Then the following are equivalent:
\begin{itemize}
    \item [(a)] $N$ satisfies the geometric criterion with respect to some complete cut system $\{D'_i\}$.
    \item[(b)] $S_{\rm min}$ satisfies $(A)$.
\end{itemize}
\end{theorem}

\begin{example} \label{ex:whitehead_reduction} Figure \ref{fig:ejemplo_complicado} showed a handlebody $N$ whose $t$--curves spell the words $S = \{x_1 x_2 x_1 x_2 x_2, x_2^{-1}x_2^{-1}x_1^{-1},x_1^{-1}x_2^{-1}\}$ with respect to the cut system $\{D_1,D_2\}$. We apply Whitehead reduction as follows. The substitution $x_1 \mapsto x_1x_2^{-1}, x_2 \mapsto x_2$ (i.e. multiplier $a = x_2^{-1}$) produces the set $S' = \{x_1x_1x_2,x_1^{-1}x_2^{-1},x_1^{-1}\}$ which has a smaller algebraic length than $S$ ($6$ versus $10$). We now substitute $x_1 \mapsto x_1, x_2 \mapsto x_1^{-1}x_2$ (multiplier $a = x_1$) and get $S'' = \{x_1x_2,x_1^{-1},x_2^{-1}\}$. It is easy to see that no substitution can reduce its length any more, so this is a minimal set $S_{\rm min}$ equivalent to $S$.
\end{example}

\begin{proof}[Proof of Theorem \ref{thm:main3_restated}] (a) $\Rightarrow$ (b) We have that each $\partial D'_i$ intersects $\bigcup \tau$ exactly twice and transversally. Thus as we travel once along $\partial D'_i$ we run into $\bigcup \tau$ exactly twice, crossing from the white region to the gray one (say) and then back to the white one. Therefore $x_i$ and $x_i^{-1}$ appear exactly once each among the unreduced words read off by the $t$--curves. (They appear with opposite exponents due to the orientation of the $t$--curves as the boundary of the gray region). Hence the set $S'$ of words read by the $t$--curves with respect to $\{D'_i\}$ satisfies condition $(A)$. We have $S' \sim S$ by Lemma \ref{lem:changebasis}, and then by Lemma \ref{lem:Qinvariant_min} the set $S_{\rm min}$ also satisfies $(A)$.

(b) $\Rightarrow$ (a) Consider all complete cut systems for $N$ which intersect $\bigcup \tau$ transversally. Since these intersections consist of finitely many points, there exists one such system $\{D'_i\}$ which minimizes the cardinality of $\left( \bigcup \partial D'_i \right) \cap \left( \bigcup \tau \right)$. Let $S'$ be the set of words read out by the $t$--curves with respect to this system. By a theorem of Zieschang (\cite[Theorem 1, p. 128]{zieschang1}, but see also \cite[p. 318]{waldhausen2} or \cite{kaneto1} for simpler proofs) all words in $S'$ are cyclically reduced and $S'$ is a minimal form of $S$. Since by assumption there is a minimal form of $S$ which satisfies $(A)$, so does $S'$ by Lemma \ref{lem:Qinvariant_min}. This directly implies that each of the curves $\partial D'_i$ either intersects $\bigcup \tau$ exactly twice (if $x_i^{\pm 1}$ appears in $S'$) or does not intersect it at all (if it does not). Notice that we are using here that the words read off by the $t$--curves in $\{D'_i\}$ are already cyclically reduced, and so they accurately record all intersections with $\bigcup D'_i$. (In general there might exist consecutive pairs of intersections with the same disk, in the form $x_i x_i^{-1}$, which would cancel out during cyclic reduction. In that case the resulting set $S'$ would underestimate the number of intersections with $\bigcup D'_i$).

We are almost finished. We only need to fix the situation when some of the $\partial D'_i$ do no intersect $\bigcup \tau$ at all. Say, for definiteness, that one of them is $\partial D'_1$. By definition $\bigcup \partial D'_i$ does not separate $\partial N$ and by assumption $\bigcup \tau$ is nonempty, so there is an arc that joins a point in $\partial D'_1$ to a point in $\bigcup \tau$ while being disjoint from the remaining $\partial D'_i$. Travel along this arc starting from its endpoint in $\partial D'_1$ until first hitting $\bigcup \tau$. Discarding the rest of the arc we have an arc $a$ whose interior is disjoint from all the $\partial D'_i$ and all the $t$--curves and whose endpoints lie one on $\partial D'_1$ and another on some $t$--curve (Figure \ref{fig:fix_intersection}.(a)).

\begin{figure}[h!]
\null\hfill
\subfigure[]{
\begin{pspicture}(0,0)(4.8,3)
	%\psgrid(0,0)(4.8,3)
	\rput[bl](0,0){\scalebox{0.7}{\includegraphics{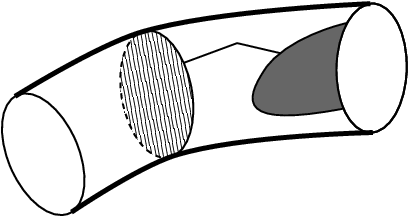}}}
	\rput[b](2.1,0.2){$D'_1$} \rput[b](2.8,2.15){$a$}
\end{pspicture}}
\hfill
\subfigure[]{
\begin{pspicture}(0,0)(4.8,3)
	%\psgrid(0,0)(4.8,3)
	\rput[bl](0,0){\scalebox{0.7}{\includegraphics{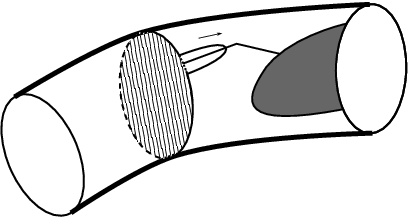}}}
\end{pspicture}}
\hfill
\subfigure[]{
\begin{pspicture}(0,0)(4.8,3)
	%\psgrid(0,0)(4.8,3)
	\rput[bl](0,0){\scalebox{0.7}{\includegraphics{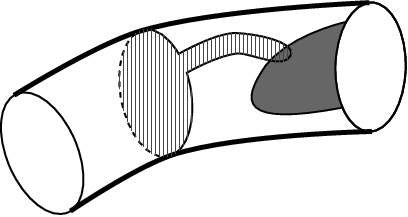}}}
\rput[b](2.1,0.2){$D''_1$}
\end{pspicture}}
\hfill\null
\caption{\label{fig:fix_intersection}}
\end{figure}

Drag a tiny portion of $\partial D'_1$ along the arc $a$ until it intersects $\tau$ in two very closely spaced points (panels (b), (c) in Figure \ref{fig:fix_intersection}). This is performed by an isotopy of $\partial N$ supported in a small neighbourhood of $a$ and can be extended to an isotopy of $N$ again supported in a small neighbourhood of $a$ (in $N$) which in particular we can take to be disjoint from all the remaining $D'_i$. The extended isotopy produces carries $D'_1$ onto a new disk $D''_1$ which together with the remaining $D'_i$ is a new cut system for $N$ where now $\partial D''_1$ intersects the system of $t$--curves transversally exactly twice. Repeating this procedure if necessary we can then obtain a new complete cut system that now satisfies the geometric criterion.
\end{proof}

We conclude with an observation: only the essential $t$--curves in $N$ are relevant as far as the geometric criterion is concerned. Ultimately this owes to the fact that any inessential $t$--curve in a coloured handlebody $N$ can be pulled off (via an isotopy) all the disks in a cut system $\{D_i\}$ without introducing any new intersections. This can be proved directly by a somewhat tedious cut-and-paste argument, but it also follows directly from the theorem above. Indeed, the set of words $S_{\rm min}$ is just the expression of (the conjugacy classes) of the $t$--curves of $N$ in a certain basis of the fundamental group $\pi_1(N)$. The inessential (i.e. contractible) $t$--curves are precisely those which correspond to the empty word $1$. Clearly empty words can be disregarded when checking the algebraic condition $(A)$, and the claim follows.

\section{Models that satisfy the geometric criterion} \label{sec:models} 

This last and brief section is motivated by the following ill-posed, but intuitively appealing, question. Suppose one observes an isolating block $N$ which is a handlebody with a complicated pattern of tangency curves. Is it likely that its maximal invariant subset has a nontrivial one-dimensional cohomology? Roughly, we need to count how many colourings of a given handlebody satisfy the geometric criterion. Of course, there are infinitely many that do (for example, any colouring with all $t$--curves inessential) and infinitely many that do not (any colouring with three or more parallel, essential $t$--curves). In spite of this, the very simple discussion in this section shows that there is a certain finiteness property in the collection of all colourings that satisfy the geometric criterion of Theorem \ref{thm:main1}.

\subsection{} Recall the handlebody $H_g$ with its standard cutting system $\{D_i\}$ shown in Figure \ref{fig:standard_hdbdy}. Suppose that it is endowed with a colouring that satisfies the geometric criterion with respect to $\{D_i\}$. Focus only on the essential $t$--curves $\{\tau'\}$, which must all intersect $\bigcup D_i$. Cutting $H_g$ open along the $D_i$ we obtain a $3$--ball $B$ with $2g$ distinguished disks $D_i^{\pm}$ on its boundary $\partial B$, each with at most two marked points in $\partial D_i^{\pm}$ corresponding to the intersection of $\partial D_i$ with $\bigcup \tau'$. The $t$--curves $\tau'$ give rise to a collection of at most $2g$ pairwise disjoint arcs in $\overline{\partial B \setminus \bigcup D_i^{\pm}}$ whose endpoints are the marked points in the boundaries of the $D_i^{\pm}$. It is not difficult to see that there are only finitely many of these collections $\mathcal{A}$ of arcs up to isotopies fixing the disks $D_i^{\pm}$. The reason is that the isotopy class of $\mathcal{A}$ depends only on how the arcs pair the marked points and what distinguished disks are enclosed by cycles of arcs. Returning to the handlebody $H_g$ by pasting back together each pair $D_i^{\pm}$ we conclude that any colouring of $H_g$ which satisfies the geometric criterion with respect to the system $\{D_i\}$ must be, after removal of its inessential $t$--curves, isotopic to one of finitely many models $\{\mathcal{M}_i^g\}$ on $H_g$.

\begin{figure}[h]
\scalebox{0.9}{\begin{pspicture}(0,0)(9.5,2.5)
	%\psgrid(0,0)(9.5,2.5)
	\rput[bl](0,0){\scalebox{0.5}{\includegraphics{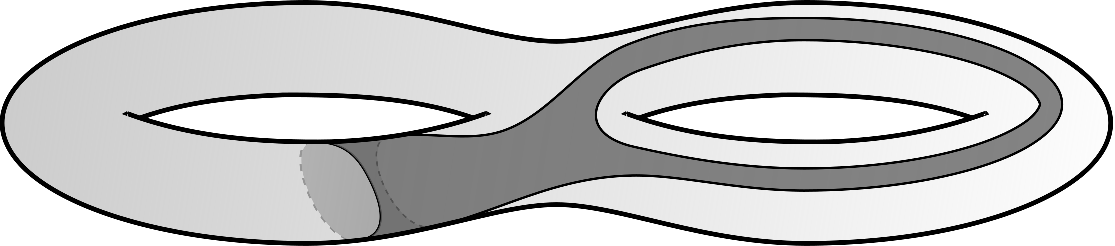}}}
\end{pspicture}}
\caption{\label{fig:example_genus2}}
\end{figure}

Figures \ref{fig:genus1} and \ref{fig:genus2} below show all possible models for genera $g = 1$ and $g = 2$. In Figure \ref{fig:genus1}.(a) the model contains no curves; there is a similar empty model for every genus (not shown in Figure \ref{fig:genus2}). Empty models correspond to the case when all the $t$--curves in the colouring of $H_g$ are inessential. We remark that, because the models $\{\mathcal{M}_i^g\}$ only show the essential $t$--curves of a colouring, some of them may not be colourable ``as is''. For instance, the model in Figure \ref{fig:genus2}.(b) needs an inessential $t$--curve to be added so that it can be coloured (see Figure \ref{fig:example_genus2}).

\begin{figure}[h]
\null\hfill
\subfigure[]{
\begin{pspicture}(0,0)(2.2,2.2)
	%\psgrid(0,0)(2.2,2.2)
	\rput[bl](0,0){\scalebox{0.2}{\includegraphics{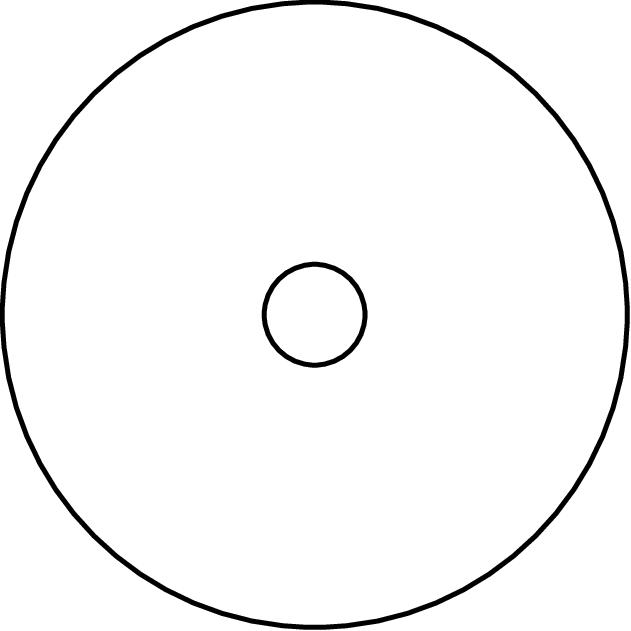}}}
\end{pspicture}}
\hfill
\subfigure[]{
\begin{pspicture}(0,0)(2.2,2.2)
	%\psgrid(0,0)(2.2,2.2)
	\rput[bl](0,0){\scalebox{0.2}{\includegraphics{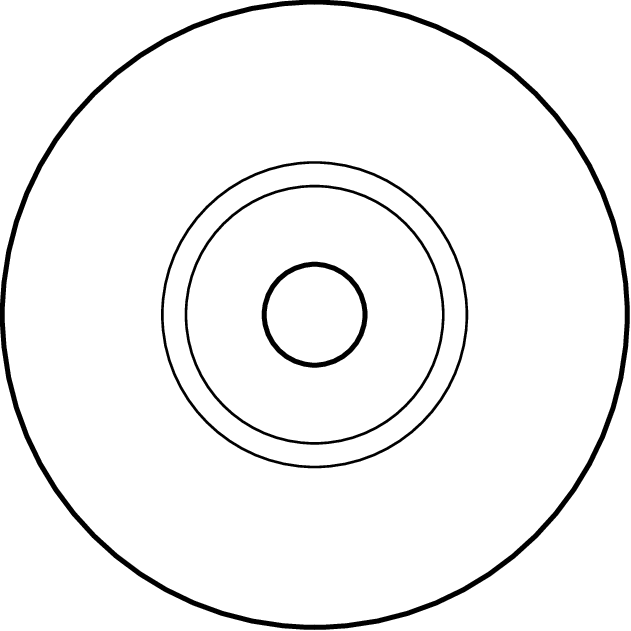}}}
\end{pspicture}}
\hfill\null
\caption{Models for genus $1$ \label{fig:genus1}}
\end{figure}

\begin{figure}[h]
\null\hfill
\subfigure[]{
\begin{pspicture}(0,0)(3.8,2.2)
	%\psgrid(0,0)(3.8,2.2)
	\rput[bl](0,0){\scalebox{0.2}{\includegraphics{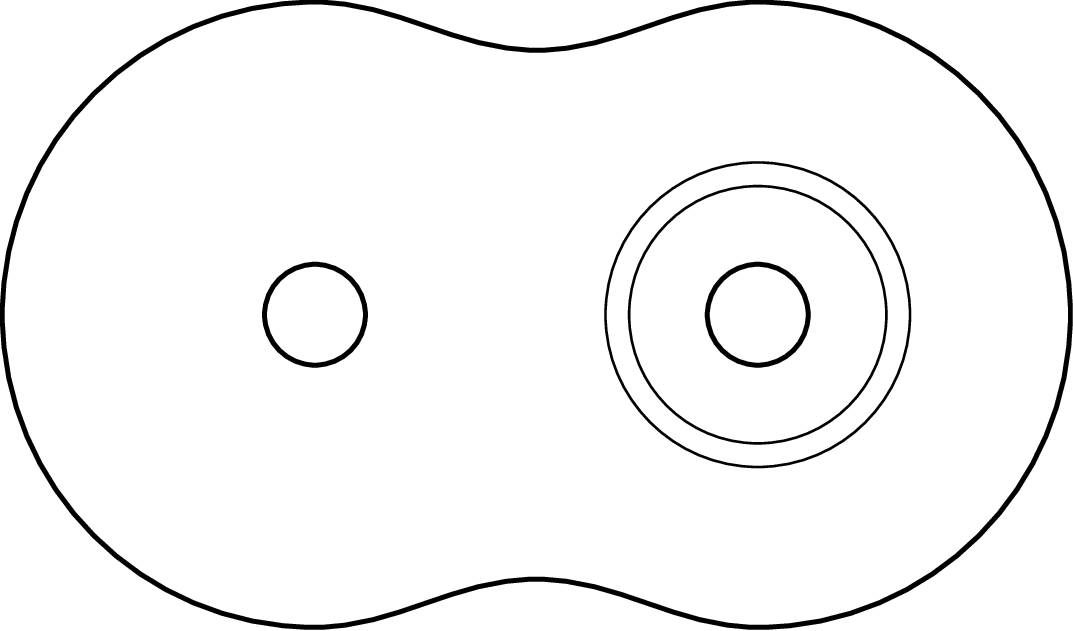}}}
\end{pspicture}}
\hfill
\subfigure[]{
\begin{pspicture}(0,0)(3.8,2.2)
	%\psgrid(0,0)(3.8,2.2)
	\rput[bl](0,0){\scalebox{0.2}{\includegraphics{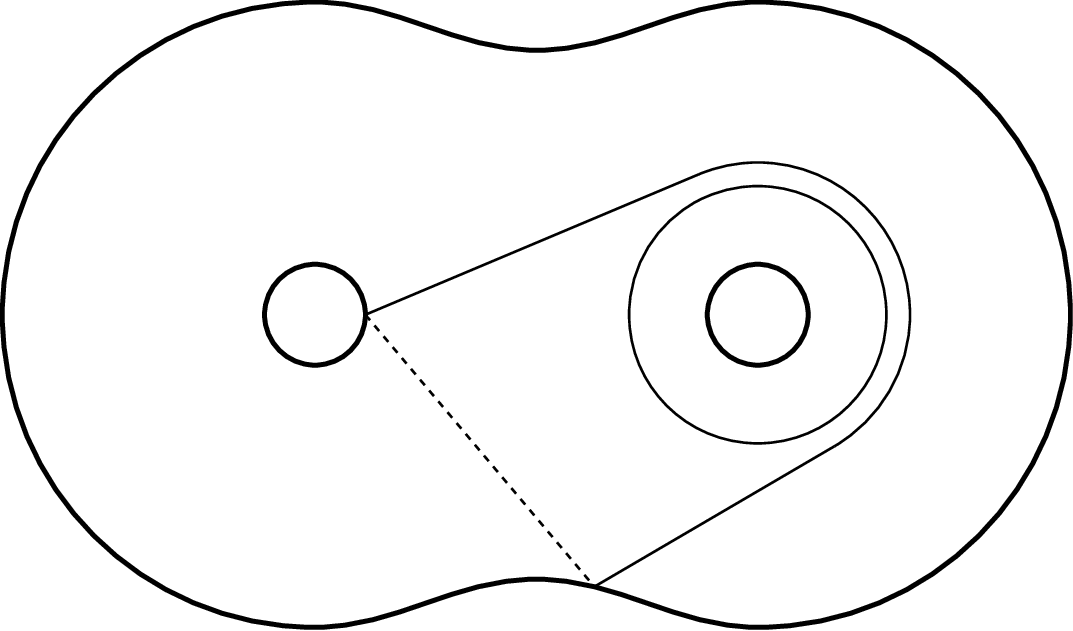}}}
\end{pspicture}}
\hfill
\subfigure[]{
\begin{pspicture}(0,0)(3.8,2.2)
	%\psgrid(0,0)(3.8,2.2)
	\rput[bl](0,0){\scalebox{0.2}{\includegraphics{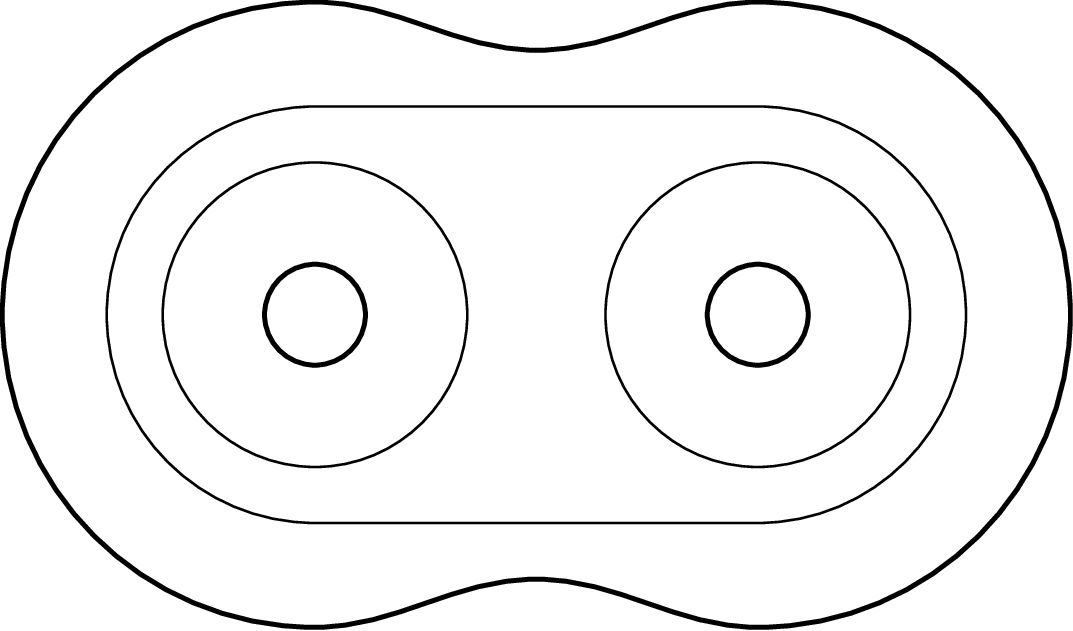}}}
\end{pspicture}}
\hfill
\subfigure[]{
\begin{pspicture}(0,0)(3.8,2.2)
	%\psgrid(0,0)(3.8,2.2)
	\rput[bl](0,0){\scalebox{0.2}{\includegraphics{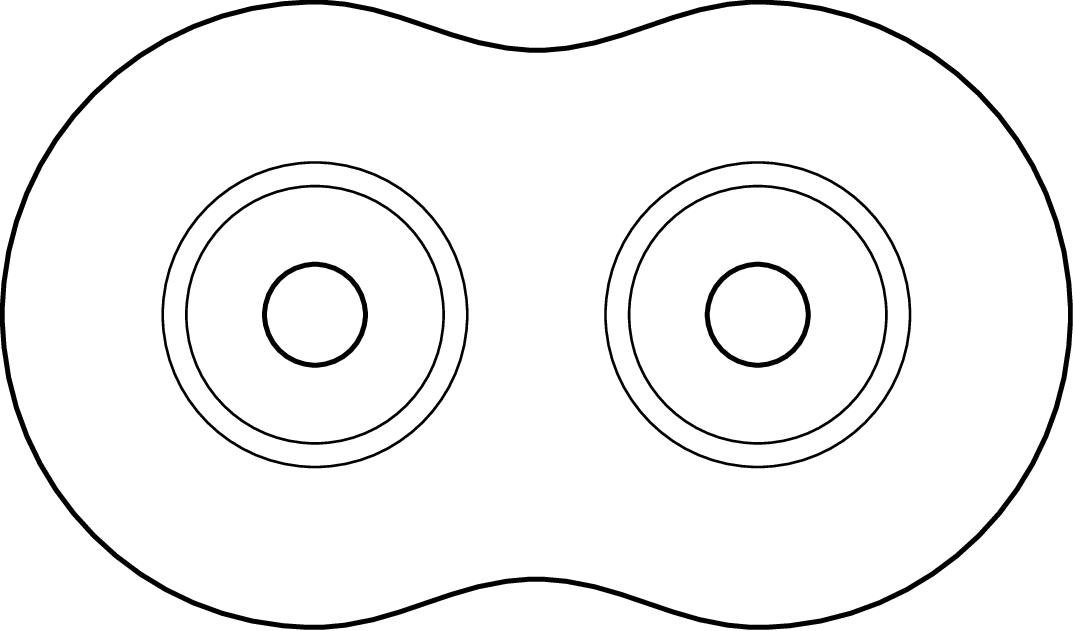}}}
\end{pspicture}}
\hfill\null
\caption{Models for genus $2$ \label{fig:genus2}}
\end{figure}

To generate these models we have essentially reversed the argument given in the previous paragraph, starting with four disjoint disks in the boundary of a $3$--ball with at most two marked points on the boundary of each and listing all possible collections of arcs joining those up to isotopy relative to the disks. One should bear in mind that this procedure is very inefficient and many of the configurations it produces have to be discarded because they include inessential curves, or they are equivalent up to a homeomorphism, or they cannot be completed (by adding inessential $t$--curves) to a collection that can be coloured.

\subsection{} We return to dynamics with the following proposition, which should be interpreted as a sort of (very weak) Hartman-Grobman theorem. For an isolated invariant set with a trivial one-dimensional cohomology it provides, up to a local topological conjugacy, ``essentially finitely many models'' for its possible isolating blocks and their entry and exit patterns. It does not say anything about the dynamics inside the block.

\begin{proposition} Suppose $N \subseteq \mathbb{R}^3$ is an isolating block with a connected boundary. Assume that its maximal invariant subset $K$ has a trivial one-dimensional cohomology. As we know by Theorem \ref{thm:main1}, $N$ must be a handlebody of some genus $g$. There exists a flow $\psi$ on $\mathbb{R}^3$ which realizes $H_g$ as an isolating block and such that:
\begin{itemize}
    \item[(i)] There exist open neighbourhoods $U$ of $N$ and $V$ of $H_g$ and a homeomorphism $h : U \longrightarrow V$ which carries $N$ onto $H_g$ and conjugates $\varphi$ and $\psi$.
    \item[(ii)] The essential tangency curves in $H_g$ accord to one of the finitely models $\{\mathcal{M}_i^g\}$ corresponding to genus $g$.
\end{itemize}
\end{proposition}
\begin{proof} By Theorem \ref{thm:main1} the system of tangency curves of $N$ satisfies the geometric condition for some complete cut system $\{D'_i\}$. There is a homeomorphism $h_1 : N \longrightarrow H_g$ which carries each $D'_i$ onto the corresponding $D_i$: this is a standard, purely topological fact about handlebodies and it is easily proved by first cutting both $N$ and $H_g$ open into $3$--balls, each with $2g$ marked disks, and constructing a homeomorphism between these two balls that matches corresponding disks. 

Copying the colouring of $N$ onto $H_g$ via $h_1$ we obtain a colouring of $H_g$ which evidently satisfies the geometric criterion with respect to $\{D_i\}$, and so there exists a homeomorphim $h_2$ of $H_g$ which carries the essential $t$--curves of the colouring onto one of the finitely many models $\mathcal{M}_i^g$. The composition $h := h_2 \circ h_1 : N \longrightarrow H_g$ then carries the essential tangency curves of $N$ onto $\mathcal{M}_i^g$. Since $N$ is tame by definition and $H_g$ is obviously tame, both can be slightly enlarged within $\mathbb{R}^3$ by adding external collars to their boundaries, obtaining neighbourhoods $U = N \cup \partial N \times [0,1)$ and $B = H_g \cup \partial H_g \times [0,1)$ of $N$ and $H_g$. $h$ can be extended to a homeomorphism $h : U \longrightarrow V$ in the obvious way. Copying the (now local) flow in $U$ to $V$ via $h$ produces a local flow $\psi$ there which can easily be extended to a complete flow in $\mathbb{R}^3$ by slowing it down to a halt near ${\rm fr}\ V$ and letting all points outside $V$ be stationary.
\end{proof}

\bibliographystyle{plain}
\bibliography{main}

\end{document}